\documentclass{article}%
\usepackage{amsfonts}
\usepackage{graphicx}
\usepackage{amsmath}
\usepackage{amssymb}%
\providecommand{\U}[1]{\protect\rule{.1in}{.1in}}
\newtheorem{theorem}{Theorem}

\newtheorem{corollary}[theorem]{Corollary}

\newtheorem{example}[theorem]{Example}

\newtheorem{lemma}[theorem]{Lemma}

\newtheorem{remark}[theorem]{Remark}

\newenvironment{proof}[1][Proof]{\textbf{#1.} }{\ \rule{0.5em}{0.5em}}
\begin{document}

\title{Creating Domain Mappings}
\author{Kendall Atkinson\\Departments of Mathematics \& Computer Science \\The University of Iowa
\and Olaf Hansen\\Department of Mathematics \\California State University San Marcos}
\maketitle

\begin{abstract}
Consider being given a mapping $\varphi:S^{d-1}%
\overset{1-1}{\underset{onto}{\longrightarrow}}\partial\Omega$, with
$\partial\Omega$ the $\left(  d-1\right)  $-dimensional smooth boundary
surface for a bounded open simply-connected region $\Omega$ in $\mathbb{R}%
^{d}$, $d\geq2$. We consider the problem of constructing an extension
$\Phi:\overline{B}_{d}\overset{1-1}{\underset{onto}{\longrightarrow}}%
\overline{\Omega}$ with $B_{d}$ the open unit ball in $\mathbb{R}^{d}$. The
mapping is also required to be continuously differentiable with a non-singular
Jacobian matrix at all points. We discuss ways of obtaining initial guesses
for such a mapping $\Phi$ and of then improving it by an iteration method.

\end{abstract}

\section{Introduction}

Consider the following problem. We are given
\begin{equation}
\varphi:\partial B_{d}\underset{onto}{\overset{1-1}{\longrightarrow}}%
\partial\Omega\label{eq1}%
\end{equation}
Notationally, $d\geq2$, $B_{d}$ is the open unit ball in $\mathbb{R}^{d}$ with
boundary $S^{d-1}=\partial B_{d}$, and $\Omega$ is an open, bounded,
simply-connected region in $\mathbb{R}^{d}$. We want to construct a
continuously differentiable extension
\begin{equation}
\Phi:\overline{B}_{d}\underset{onto}{\overset{1-1}{\longrightarrow}}%
\overline{\Omega}\label{eq2}%
\end{equation}
such that%
\begin{align}
\left.  \Phi\right\vert _{S^{d-1}}  &  =\varphi\smallskip\label{eq3}\\
\det\left(  D\Phi\left(  x\right)  \right)   &  \neq0,\quad\quad x\in
\overline{B}_{d}\label{eq4}%
\end{align}
$D\Phi$ denotes the $d\times d$ Jacobian of $\Phi$,%
\[
\left(  D\Phi\left(  x\right)  \right)  _{i,j}=\frac{\partial\Phi_{i}%
(x)}{\partial x_{j}},\quad\quad x\in\overline{B}_{d}%
\]
\ As a particular case, let $d=2$ and consider extending a smooth mapping%
\[
\varphi:S^{1}\underset{onto}{\overset{1-1}{\longrightarrow}}\partial\Omega
\]
with $\Omega$ an open, bounded region in $\mathbb{R}^{2}$ and $\varphi$ a
smooth mappinng. In this case, a conformal mapping will give a desirable
solution $\Phi$; but finding the conformal mapping is often nontrivial.. In
addition, our eventual applications need the Jacobian $D\Phi$ (see
\cite{ACH2010}, \cite{AH2010}, \cite{AHC2011}), and obtaining $D\Phi$ is
difficult with most methods for constructing conformal mappings. As an
example, let $\varphi$ define an ellipse,%
\[
\varphi\left(  \cos\theta,\sin\theta\right)  =\left(  a\cos\theta,b\sin
\theta\right)  ,\quad0\leq\theta\leq2\pi
\]
with $a,b>0$. The conformal mapping has a complicated construction requiring
elliptic functions (e.g. see \cite[\S 5]{atkinson-han2004}), whereas the much
simpler mapping%
\[
\Phi\left(  x_{1},x_{2}\right)  =\left(  ax_{1},bx_{2}\right)  ,\quad
x\in\overline{B}_{2}%
\]
is sufficient for most applications. Also, for $d>2$, constructing a conformal
mapping is no longer an option.

In \S 2 we consider various methods that can be used to construct $\Phi$, with
much of our work considering regions $\Omega$ that are `star-like' with
respect to the origin:%
\begin{align}
\varphi\left(  x\right)   &  =\rho\left(  x\right)  x,\quad\quad x\in
S^{d-1}\smallskip\label{eq5}\\
\rho &  :S^{d-1}\underset{onto}{\overset{1-1}{\longrightarrow}}\mathbb{R}%
_{>0}\nonumber
\end{align}
For convex regions $\Omega$, an integration based formula is given, analyzed,
and illustrated in \S 3. In \S 4 we present an optimization based iteration
method for improving `initial guesses' for $\Phi$. Most of the presentation
will be for the planar case ($d=2$); the case of $d=3$ is presented in
\S \ref{3D}.

\section{Constructions of $\Phi$\label{constructions}}

Let $\Omega$ be star-like with respect to the origin. We begin with an
illustration of an apparently simple construction that does not work in most
cases. Consider that our initial mapping $\varphi$ is of the form (\ref{eq5}).
Define%
\begin{align}
\Phi\left(  x,y\right)   &  =r\widehat{\rho}\left(  \theta\right)  \left(
\cos\theta,\sin\theta\right)  ,\quad\quad0\leq r\leq1\label{eq6}\\
&  =\widehat{\rho}\left(  \theta\right)  x\label{eq7}%
\end{align}
with $x=\left(  r\cos\theta,r\sin\theta\right)  $ and $\widehat{\rho}\left(
\theta\right)  =\rho\left(  x\right)  $ a periodic nonzero positive function
over $\left[  0,2\pi\right]  $. This mapping $\Phi$ has differentiability
problems at the origin $\left(  0,0\right)  $. To see this, we need to find
the derivatives of $\widehat{\rho}(\theta)$ with respect to $x$ and $y$. Use%
\[
\theta=\tan^{-1}\left(  y/x\right)  ,\quad\quad y>0,\quad x\neq0
\]
and an appropriate modification for points $\left(  x,y\right)  $ in the lower
half-plane. We find the derivatives of $\widehat{\rho}$ using%
\[
\frac{\partial\widehat{\rho}(\theta)}{\partial x}=\widehat{\rho}^{\prime
}(\theta)\frac{\partial\theta}{\partial x},\quad\quad\frac{\partial
\widehat{\rho}(\theta)}{\partial y}=\widehat{\rho}^{\prime}(\theta
)\frac{\partial\theta}{\partial y}%
\]
Then%
\[
\frac{\partial\theta}{\partial x}=\frac{-y}{x^{2}+y^{2}},\quad\quad
\frac{\partial\theta}{\partial y}=\frac{x}{x^{2}+y^{2}}%
\]

Using these,
\begin{align*}
\frac{\partial\Phi}{\partial x}  &  =\left(  \widehat{\rho}(\theta)-\frac
{xy}{x^{2}+y^{2}}\widehat{\rho}^{\prime}(\theta),\frac{-y^{2}}{x^{2}+y^{2}%
}\widehat{\rho}^{\prime}(\theta)\right)  \medskip\\
\frac{\partial\Phi}{\partial y}  &  =\left(  \frac{x^{2}}{x^{2}+y^{2}%
}\widehat{\rho}^{\prime}(\theta),\rho(\theta)+\frac{xy}{x^{2}+y^{2}%
}\widehat{\rho}^{\prime}(\theta)\right)
\end{align*}
The functions
\[
\frac{x^{2}}{x^{2}+y^{2}},\quad\frac{xy}{x^{2}+y^{2}},\quad\frac{y^{2}}%
{x^{2}+y^{2}}%
\]
are not continuous at the origin. This concludes our demonstration that the
extension $\Phi$ of (\ref{eq6}) does not work.

\subsection{Harmonic mappings\label{HarmMap}}

As our first construction method for $\Phi$, consider the more general problem
of extending to all of $B_{d}$ a real or complex valued function $f$ defined
on the boundary of $B_{d}$. Expand $f$ in a Fourier series,%
\begin{equation}
f(\theta)=\frac{1}{2}a_{0}+\sum_{n=1}^{\infty}a_{n}\cos\left(  n\theta\right)
+b_{n}\sin\left(  n\theta\right) \label{e4}%
\end{equation}
Define $F$ on $B$ using%
\begin{equation}
F\left(  x\right)  =\frac{1}{2}a_{0}+\sum_{n=1}^{\infty}r^{n}\left[  a_{n}%
\cos\left(  n\theta\right)  +b_{n}\sin\left(  n\theta\right)  \right]
\label{e5}%
\end{equation}
with $x=\left(  r\cos\theta,r\sin\theta\right)  $. Note that this is the
solution to the Dirichlet problem for Laplace's equation on the unit disk,
with the boundary data given by $f(\theta)$, $0\leq\theta\leq2\pi$.%

\begin{figure}[h]%
\centering
\includegraphics[
height=3in,
width=3.9998in
]%
{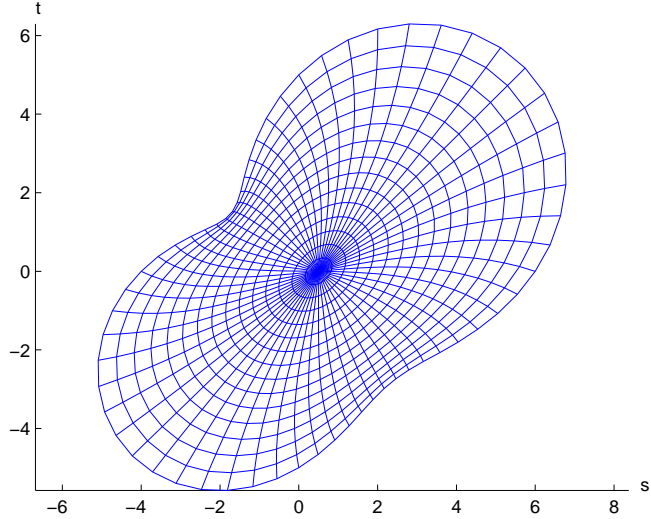}%
\caption{Starlike region with $\protect\widehat{\rho}$ defined by (\ref{e20})
with $a=5$}%
\label{fig1}%
\end{figure}

It is straightforward to show that $F$\ is infinitely differentiable for
$\left\vert x\right\vert <1$, a well-known result. In particular,
\begin{equation}
\frac{\partial F}{\partial x_{1}}=a_{1}+\sum_{m=1}^{\infty}\left(  m+1\right)
r^{m}\left[  a_{m+1}\cos m\theta+b_{m+1}\sin m\theta\right] \label{e9}%
\end{equation}%
\begin{equation}
\frac{\partial F}{\partial x_{2}}=b_{1}+\sum_{m=1}^{\infty}\left(  m+1\right)
r^{m}\left[  -a_{m+1}\sin m\theta+b_{m+1}\cos m\theta\right] \label{e10}%
\end{equation}
Depending on the speed of convergence of (\ref{e4}), we have the partial
derivatives of $F(x)$ are continuous over $\overline{B}_{2}$. In particular,
if we have
\[
\sum_{n=1}^{\infty}n\left\vert a_{n}\right\vert <\infty,\quad\quad\sum
_{n=1}^{\infty}n\left\vert b_{n}\right\vert <\infty
\]
then $\partial F/\partial x_{1}$ and $\partial F/\partial x_{2}$ are
continuous over $\overline{B}_{2}$.

Given a boundary function%
\begin{equation}
\varphi(\theta)=\left(  \varphi_{1}\left(  \theta\right)  ,\varphi_{2}\left(
\theta\right)  \right)  ,\quad\quad0\leq\theta\leq2\pi,\label{e11}%
\end{equation}
we can expand each component to all of $B_{2}$ using the above construction in
(\ref{e5}), obtaining a function $\Phi$ defined on $B_{d}$ into $\mathbb{R}%
^{2}$. A similar construction can be used for higher dimensions using an
expansion with spherical harmonics. It is unknown whether the mapping $\Phi$
obtained in this way is a one-to-one mapping from $B_{2}$ onto $\Omega$, even
if $\Omega$ is convex.

The method can be implemented as follows.

\begin{itemize}
\item Truncate the Fourier series for each of the functions $\varphi
_{k}\left(  \theta\right)  $, $k=1,2,$ say to trigonometric polynomials of
degree $n$.

\item Approximate the Fourier coefficients $\left\{  a_{j}\right\}  $ and
$\left\{  b_{j}\right\}  $ for the truncated series.

\item Define the extensions $\Phi_{k}\left(  x\right)  $ in analogy with
(\ref{e5}).
\end{itemize}%

\begin{figure}[h]%
\centering
\includegraphics[
height=3in,
width=3.9998in
]%
{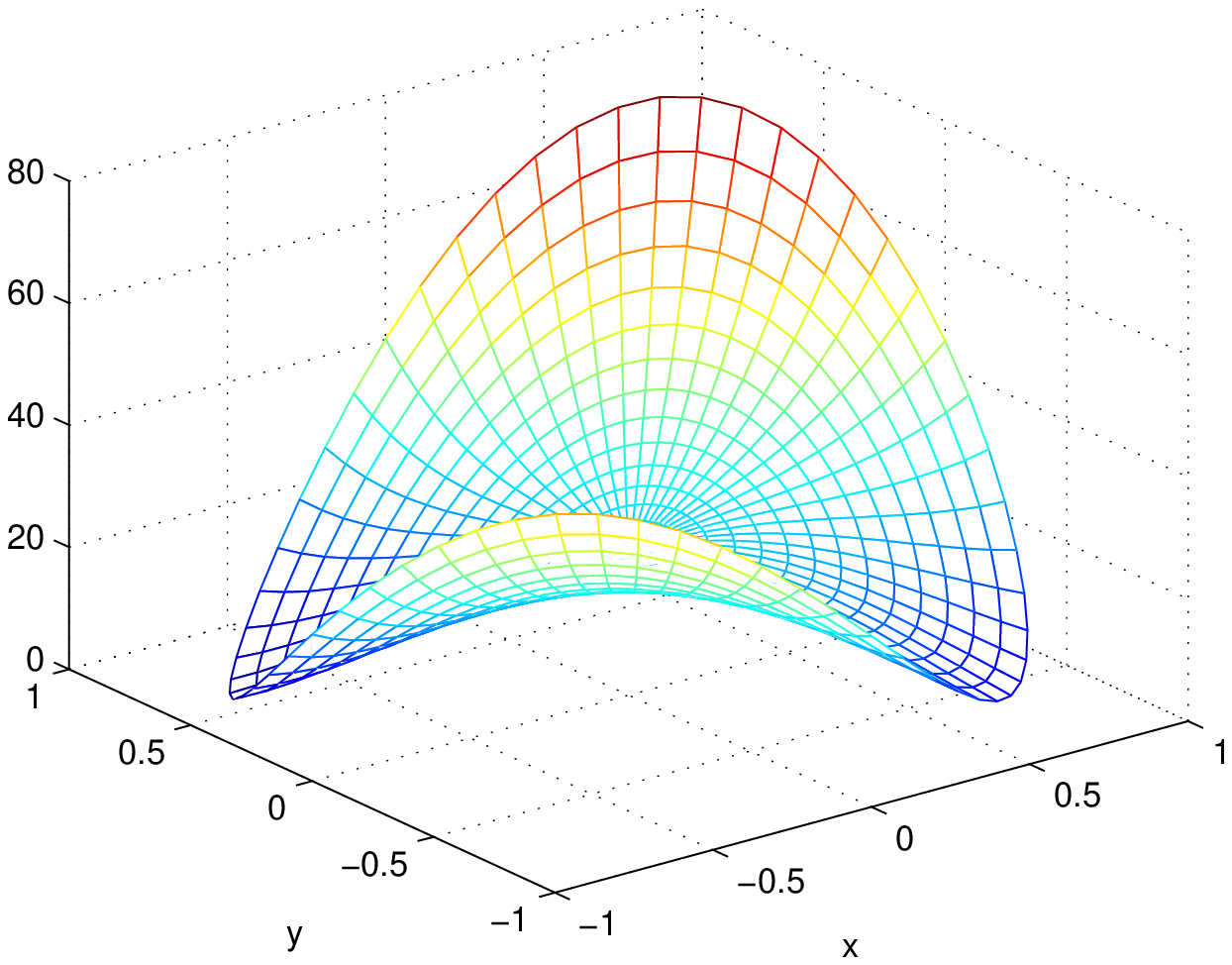}%
\caption{The Jacobian for (\ref{e20}) with $0.905\leq\left\vert \det\left(
D\Phi\left(  x\right)  \right)  \right\vert \leq75.314$}%
\label{fig2}%
\end{figure}

\begin{example}
\label{Exam1}Choose%
\begin{equation}
\rho(\theta)=a+\cos\theta+2\sin2\theta\label{e20}%
\end{equation}
with $a$ chosen greater than the maximum of $\left\vert \cos\theta
+2\sin2\theta\right\vert $ for $0\leq\theta\leq2\pi$, approximately 2.2361.
\ Note that $\rho(\theta)\cos\theta$ and $\rho(\theta)\sin\theta$ are trig
polynomials of degree 3. Begin by choosing $a=5$. Choosing \ $n=3$, we obtain
the graphs in Figures \ref{fig1} and \ref{fig2}. Figure \ref{fig1}
demonstrates the mapping by showing the images in $\Omega$ of the circles
$r=j/p$, $j=0,,\dots,p$ and the azimuthal lines $\theta=\pi j/p$,
$j=1,\dots,2p$, $p=15$. For the numerical evaluation of the Fourier
coefficients, the trapezoidal rule with $10$ nodes was used. Figure \ref{fig2}
shows $\left\vert \det\left(  D\Phi\left(  x\right)  \right)  \right\vert $
The figures illustrate that this $\Phi$ is a satisfactory mapping. However, it
is possible to improve on this mapping in the sense of reducing the ratio%
\begin{equation}
\Lambda\left(  \Phi\right)  \equiv\frac{\max\limits_{x\in\overline{B}_{2}%
}\left\vert \det\left(  D\Phi\left(  x\right)  \right)  \right\vert }%
{\min\limits_{x\in\overline{B}_{2}}\left\vert \det\left(  D\Phi\left(
x\right)  \right)  \right\vert }\label{eq10}%
\end{equation}
For the present case, $\Lambda\left(  \Phi\right)  =100.7$. An iteration
method for decreasing the size of $\Lambda\left(  \Phi\right)  $ is discussed
in \S 4. As a side-note, in the planar graphics throughout this paper we label
the axes over the unit disk as $x$ and $y$, and over $\Omega$, we label them
as $s$ and $t$.
\end{example}

In contrast to this example, when choosing $a=3$ in (\ref{e20}) the mapping
$\Phi$ derived in the same manner is neither one-to-one nor onto. Another
method is needed to generate a mapping $\Phi$ which satisfies (\ref{e20})
(\ref{eq2})-(\ref{eq4}).%
\begin{figure}[tb]%
\centering
\includegraphics[
height=3in,
width=3.9998in
]%
{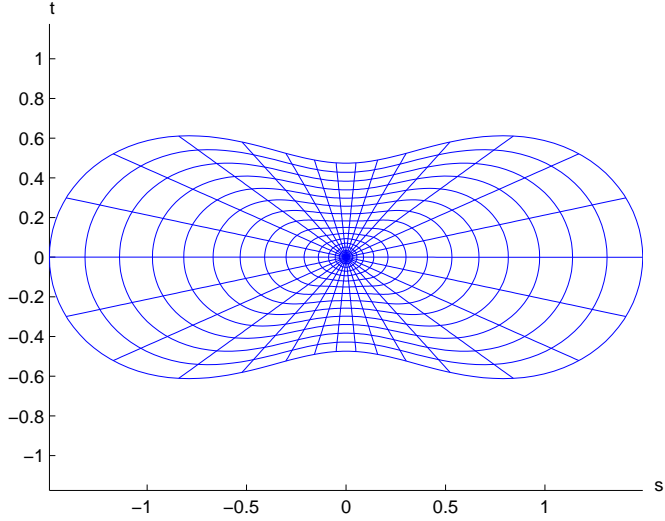}%
\caption{Starlike Cassini region with $\protect\widehat{\rho}$ defined in
(\ref{eq20}) with $a=1.5$}%
\label{ovals3}%
\end{figure}

\subsection{Using $C^{\infty}$-modification functions}

Let $\left(  x,y\right)  =\left(  r\cos\theta,r\sin\theta\right)  $, $0\leq
r\leq1$. As earlier in (\ref{eq5}), consider $\Omega$ as star-like with
respect to the origin \ Introduce the function
\begin{equation}
T(r;\kappa)=\exp\left(  \kappa\left(  1-\frac{1}{r}\right)  \right)
,\quad\quad0<r\leq1\label{e39}%
\end{equation}
with $\kappa>0$ and $T(0,\kappa)=0$. Define $\Phi$ by%
\begin{equation}
s=\Phi\left(  x;\kappa,\omega\right)  =\left[  T\left(  r;\kappa\right)
\widehat{\rho}\left(  \theta\right)  +\left(  1-T\left(  r;\kappa\right)
\right)  \omega\right]  x,\quad\quad x\in B_{2}\label{e40}%
\end{equation}
with $x=r\left(  \cos\theta,\sin\theta\right)  $, for $0\leq r\leq
1,\ 0\leq\theta\leq2\pi$, with some $\omega>0$. This is an attempt to fix the
lack of differentiability at (0,0) of the mapping (\ref{eq6})-(\ref{eq7}). As
$r$ decreases to $0$, we have $\Phi\left(  x\right)  \approx\omega x$. Thus
the Jacobian of $\Phi$ is nonzero around $\left(  0,0\right)  $. The constants
$\kappa,\omega$ are to be used as additional design parameters.

The number $\omega$ should be chosen so as to also ensure the mapping $\Phi$
is 1-1 and onto. Begin by finding a disk centered at $\left(  0,0\right)  $
that is entirely included in the open set $\Omega$, and say its radius is
$\omega_{0},$ or define%
\[
\omega_{0}=\min_{0\leq\theta\leq2\pi}\widehat{\rho}\left(  \theta\right)
\]
Then choose $\omega\in(0,\omega_{0})$. To see this is satisfactory, write%
\[%
\begin{tabular}
[c]{c}%
$\Phi\left(  x;\kappa,\omega\right)  =f(r)\left(  \cos\theta,\sin
\theta\right)  $\smallskip\\
$f(r)=r\left[  T\left(  r;\kappa\right)  \widehat{\rho}\left(  \theta\right)
+\left(  1-T\left(  r;\kappa\right)  \right)  \omega\right]  $%
\end{tabular}
\]
fixing $\theta\in\left[  0,2\pi\right]  $. Immediately, $f\left(  0\right)  =0
$, $f\left(  1\right)  =\widehat{\rho}\left(  \theta\right)  $. By a
straightforward computation,
\[
f^{\prime}\left(  r\right)  =\frac{1}{r}\left\{  T\left[  \left(
\widehat{\rho}-\omega\right)  \left(  r+1\right)  \right]  +\omega\right\}
\]
where $T=T\left(  r;\kappa\right)  $ and $\widehat{\rho}=\widehat{\rho}\left(
\theta\right)  $. The assumption $0<\omega<\omega_{0}$ then implies
\[
f^{\prime}\left(  r\right)  >0,\quad\quad0<r\leq1
\]
Thus the mapping $f:\left[  0,1\right]  \rightarrow\left[  0,\widehat{\rho
}\left(  \theta\right)  \right]  $ is 1-1 and onto, and from this
$\Phi:\overline{B}_{2}\rightarrow\overline{\Omega}$ is 1-1 and onto for the
definition in (\ref{e40}).

This definition of $\Phi$ satisfies (\ref{eq2})-(\ref{eq4}), but often leads
to a large value for the ratio $\Lambda\left(  \Phi\right)  $ of (\ref{eq10}).
It can be used as an initial choice for a $\Phi$ that can be improved by the
iteration method defined in \S 4.

\begin{example}
Consider the starlike region with
\begin{equation}
\widehat{\rho}\left(  \theta\right)  =\sqrt{\cos\left(  2\theta\right)
+\sqrt{a-\sin^{2}\left(  2\theta\right)  }}\label{eq20}%
\end{equation}
with $a>1$. The region $\Omega$ is called an `\textit{ovals of Cassini'}. We
give an example with $a=1.5$, $\left(  \kappa,\omega\right)  =\left(
1.0,0.5\right)  $. Figure \ref{ovals3} is the analog of Figure \ref{fig1}. For
the Jacobian,
\begin{align*}
\min_{r\leq1}D\Phi(x,y)  &  =0.0625\smallskip\\
\max_{r\leq1}D\Phi(x,y)  &  =4.0766
\end{align*}
The ratio $\Lambda\left(  \Phi\right)  =65.2$ is large and can be made
smaller; see Example \ref{Exam_IT4}.
\end{example}

A variation to (\ref{e40}) begins by finding a \ closed disk about the origin
that is contained wholly in the interior of $\Omega$. Say the closed disk is
of radius $\delta,$ $0<\delta<1$. Then define%
\begin{equation}%
\begin{array}
[c]{l}%
\Phi\left(  x;\kappa,\omega\right)  \smallskip\\
\quad=\left\{
\begin{array}
[c]{ll}%
x, & 0\leq r\leq\delta\medskip\\
\left[  T\left(  \dfrac{r-\delta}{1-\delta},\kappa\right)  \rho\left(
\theta\right)  +\left(  1-T\left(  \dfrac{r-\delta}{1-\delta},\kappa\right)
\right)  \right]  x, & \delta<r\leq1,
\end{array}
\right.
\end{array}
\label{e50}%
\end{equation}
where $x=r\left(  \cos\theta,\sin\theta\right)  $. Then the Jacobian $D\Phi$
around the origin is simply the identity matrix, and this ensures that $\det
D\Phi(x)\neq0$ for $x\in\overline{B}_{d}$. Experimentation is recommended on
the use of either (\ref{e40}) or (\ref{e50}), including how to choose $\kappa
$, $\omega$, and $\delta$.

The methods of this section generalize easily to the determination of an
extension $\Phi:B_{3}\underset{onto}{\overset{1-1}{\longrightarrow}}\Omega$
for the given boundary mapping%
\[
\varphi:\partial B_{3}\underset{onto}{\overset{1-1}{\longrightarrow}}%
\partial\Omega
\]
Examples of such are illustrated later in \S \ref{3D}.

\section{An integration-based mapping formula\label{integinterp}}

Begin by considering a point $P=r(\cos\alpha,\sin\alpha)\in B_{2}$,
$r\in\lbrack0,1)$, $\alpha\in\lbrack0,2\pi)$. Given an angle $\alpha\leq
\theta<\pi+\alpha$, draw a line $L$ through $P$ at an angle of $\theta$ with
respect to the positive $x_{1}$-axis. Let $P_{+}(\theta)$ and $P_{-}(\theta)$
denote the intersection of this line with the unit disk. These points will
have the form%
\begin{equation}%
\begin{array}
[c]{c}%
P_{+}\left(  \theta\right)  =P+r_{+}\left(  \theta\right)  \boldsymbol{\eta
},\\
P_{-}\left(  \theta\right)  =P-r_{-}\left(  \theta\right)  \boldsymbol{\eta}.
\end{array}
\label{e3z}%
\end{equation}
with
\begin{equation}
\boldsymbol{\eta}=\left(  \cos\theta,\sin\theta\right)  ,\quad\quad\alpha
\leq\theta<\pi+\alpha.\label{e3}%
\end{equation}
We choose $r_{+}\left(  \theta\right)  $ and $r_{-}\left(  \theta\right)  $ to
be such that%
\[
\left\vert P_{+}(\theta)\right\vert =\left\vert P+r_{+}\left(  \theta\right)
\boldsymbol{\eta}\right\vert =1,\quad\quad\left\vert P_{-}(\theta)\right\vert
=\left\vert P-r_{-}\left(  \theta\right)  \boldsymbol{\eta}\right\vert =1
\]
and%
\[
\left\vert r_{+}\left(  \theta\right)  \right\vert =\left\vert P-P_{+}%
(\theta)\right\vert ,\quad\left\vert r_{-}\left(  \theta\right)  \right\vert
=\left\vert P-P_{-}(\theta)\right\vert .
\]

Define
\begin{equation}
\varphi_{\ast}\left(  \theta\right)  =\varphi\left(  P_{+}\left(
\theta\right)  \right)  -r_{+}(\theta)\frac{\varphi\left(  P_{+}\left(
\theta\right)  \right)  -\varphi\left(  P_{-}\left(  \theta\right)  \right)
}{r_{+}(\theta)+r_{-}(\theta)}\label{eq24}%
\end{equation}
using linear interpolation along the line $L$. Here and in the following we
always identify the function $\varphi$ on the boundary of the unit disk with a
$2\pi$ periodic function on the real number line. \ Then define%
\begin{equation}
\Phi\left(  P\right)  =\frac{1}{\pi}\int_{\alpha}^{\alpha+\pi}\varphi_{\ast
}\left(  \theta\right)  \,d\theta\label{eq25}%
\end{equation}
We study the construction and properties of $\Phi$ in the following two sections.

\subsection{Constructing $\Phi$}

The most important construction is the calculation of $P_{+}\left(
\theta\right)  $\ and $P_{-}\left(  \theta\right)  $\textbf{. }We want to find
two points $\boldsymbol{\gamma}$ \ that are the intersection of $\partial
B_{2}$ and the straight line $L$ through $P$ in the direction
$\boldsymbol{\eta}$, $\left\vert \boldsymbol{\eta}\right\vert =1$. \ Since
$P\in\operatorname*{int}\left(  B_{2}\right)  $, we have $\left\vert
P\right\vert <1$. We want to find
\[
\boldsymbol{\gamma}=P+r\boldsymbol{\eta},\quad\quad\left\vert
\boldsymbol{\gamma}\right\vert =1
\]
with $\boldsymbol{\eta}$ denoting the direction from $P$ as noted earlier.
With the assumption (\ref{e3}) on \ $\boldsymbol{\eta}$, we have%
\[%
\begin{array}
[c]{cl}%
0\leq P\cdot\boldsymbol{\eta}\leq\left\vert P\right\vert , & \alpha\leq
\theta\leq\alpha+\frac{1}{2}\pi\smallskip\\
0\leq-P\cdot\boldsymbol{\eta}\leq\left\vert P\right\vert , & \alpha+\frac
{1}{2}\pi\leq\theta\leq\alpha+\pi
\end{array}
\]
Using $\boldsymbol{\gamma}\cdot\boldsymbol{\gamma}=1$,%
\[
\left\vert P+r\boldsymbol{\eta}\right\vert ^{2}=P\cdot P+2rP\cdot
\boldsymbol{\eta}+r^{2}=1
\]%
\[
r^{2}+2rP\cdot\boldsymbol{\eta}+\underset{<0}{\underbrace{P\cdot P-1}}=0
\]
which implies the roots are real and nonzero. Thus the formula%
\begin{align*}
r  &  =-\boldsymbol{\eta}\cdot P\pm\sqrt{\left(  P\cdot\boldsymbol{\eta
}\right)  ^{2}+1-P\cdot P}\\
&  =r\cos(\theta-\alpha)\pm\sqrt{1-r^{2}\sin^{2}(\theta-\alpha)}%
\end{align*}
defines two real roots. Here we see that%
\[%
\begin{tabular}
[c]{c}%
$\left\vert P\cdot\boldsymbol{\eta}\right\vert \leq\left\vert P\right\vert $\\
$\left(  P\cdot\boldsymbol{\eta}\right)  ^{2}-P\cdot P\leq0$\\
$\left(  P\cdot\boldsymbol{\eta}\right)  ^{2}+1-P\cdot P\leq1$%
\end{tabular}
\]
So
\begin{align*}
r_{-}  &  =P\cdot\boldsymbol{\eta}+\sqrt{\left(  P\cdot\boldsymbol{\eta
}\right)  ^{2}+1-P\cdot P}\\
&  =r\cos(\theta-\alpha)+\sqrt{1-r^{2}\sin^{2}(\theta-\alpha)}\\
r_{+}  &  =-P\cdot\boldsymbol{\eta}+\sqrt{\left(  P\cdot\boldsymbol{\eta
}\right)  ^{2}+1-P\cdot P}\\
&  =-r\cos(\theta-\alpha)+\sqrt{1-r^{2}\sin^{2}(\theta-\alpha)}%
\end{align*}
It is immediate that%
\begin{align*}
r_{-}+r_{+}  &  =2\sqrt{\left(  P\cdot\boldsymbol{\eta}\right)  ^{2}+1-P\cdot
P}\\
&  =2\sqrt{1-r^{2}\sin^{2}(\theta-\alpha)}%
\end{align*}
and therefore the denominator in the formula (\ref{eq24}) for $\varphi_{\ast
}\left(  \theta\right)  $ is zero if and only if \ $\left\vert P\right\vert
=1$ and $P\perp\boldsymbol{\eta}$, a case not allowed in our construction.

Using $r_{-}$ and $r_{+}$ in (\ref{e3z}), we can construct $\varphi_{\ast
}\left(  \theta\right)  $ using (\ref{eq24}), and this is then used in
obtaining the mapping $\Phi\left(  P\right)  $ of (\ref{eq25}). This formula
is approximated using numerical integration with the trapezoidal rule. We
illustrate this later in the section.%

\begin{figure}[tb]%
\centering
\includegraphics[
height=3in,
width=3.9998in
]%
{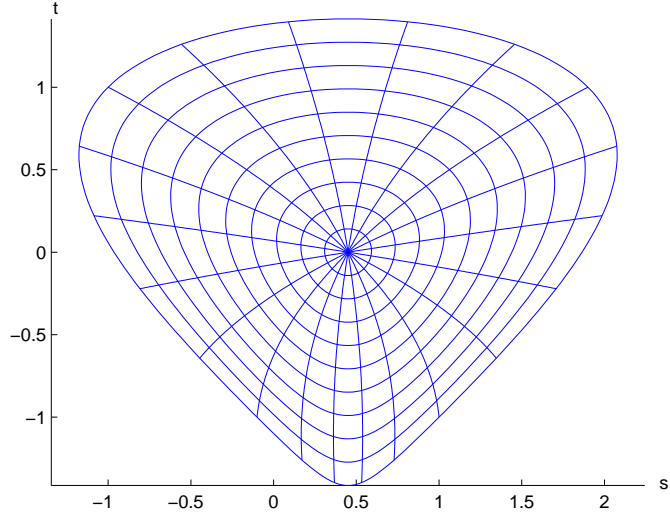}%
\caption{The mapping $\Phi$ for boundary (\ref{eq30}) with $a=0.9$}%
\label{fig5}%
\end{figure}
%

\begin{figure}[tb]%
\centering
\includegraphics[
height=3in,
width=3.9998in
]%
{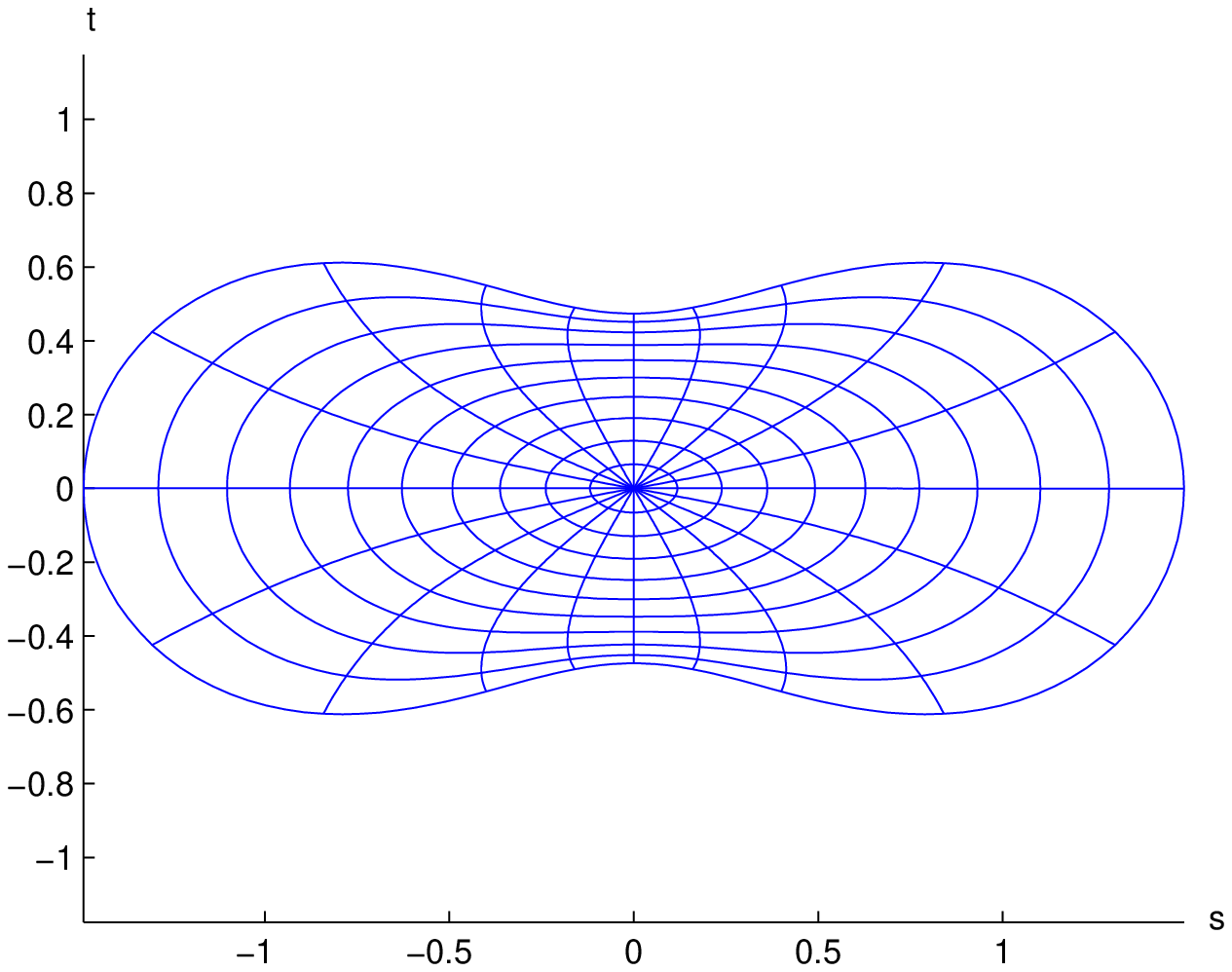}%
\caption{The mapping $\Phi$ for boundary (\ref{eq20}) with $a=1.5$}%
\label{fig6}%
\end{figure}

To further simplify the analysis of the mapping $\Phi$ of (\ref{eq25}), we
assume for a moment that $\alpha=0$, so the point $P$ is located on the
positive x--axis. Our next goal is to determine the respective angles between
$P_{+}(\theta)$ and $P_{-}(\theta)$ and the positive x-axis. We denote these
angles by $\psi^{+}$ and $\psi^{-}$, respectively. Using the law of cosines in
the triangle given by the origin, $P$, and $P^{+}$ we obtain
\begin{align}
(r_{+})^{2}  &  =r^{2}+1-2r\cos(\psi^{+})\nonumber\\
2r\cos(\psi^{+})  &  =r^{2}+1-(r_{+})^{2}\nonumber\\
&  =r^{2}+1-\left(  -r\cos(\theta)+\sqrt{1-r^{2}\sin^{2}\theta}\right)
^{2}\nonumber\\
&  =2r^{2}\sin^{2}\theta+2r\cos(\theta)\sqrt{1-r^{2}\sin^{2}\theta}\nonumber
\end{align}
which implies%
\begin{equation}
\psi_{+}=\psi_{+}(r,\theta)=\arccos\left(  r\sin^{2}\theta+\cos(\theta
)\sqrt{1-r^{2}\sin^{2}\theta}\right) \label{eq2000}%
\end{equation}
where we use the function $\arccos:[-1,1]\mapsto\lbrack0,\pi]$. Similarly we
get
\begin{equation}
\psi_{-}\equiv\psi_{-}(r,\theta)=\widetilde{\arccos}\left(  r\sin^{2}%
\theta-\cos(\theta)\sqrt{1-r^{2}\sin^{2}\theta}\right) \label{eq2001}%
\end{equation}
where we use the function $\widetilde{\arccos}:[-1,1]\mapsto\lbrack\pi,2\pi]
$,
\[
\widetilde{\arccos}(x)=2\pi-\arccos(x)
\]
Using the functions $\psi_{-}$ and $\psi_{+}$ we can rewrite $\varphi_{\ast
}(\theta)$, see (\ref{eq24}), in the following way:
\begin{align*}
\varphi_{\ast}(r,\theta)  &  =\frac{1}{2}\left(  1+\frac{r\cos(\theta)}%
{\sqrt{1-r^{2}\sin^{2}\theta}}\right)  \varphi(\psi_{+}(r,\theta))\\
&  +\frac{1}{2}\left(  1-\frac{r\cos(\theta)}{\sqrt{1-r^{2}\sin^{2}\theta}%
}\right)  \varphi(\psi_{-}(r,\theta))
\end{align*}
This allows us to write formula (\ref{eq25}) more explicitly in the following
way:
\begin{align}
\Phi(P)  &  =\frac{1}{2\pi}%
{\displaystyle\int_{0}^{\pi}}
\left(  1+\frac{r\cos(\theta)}{\sqrt{1-r^{2}\sin^{2}\theta}}\right)
\varphi(\psi_{+}(r,\theta))d\theta\nonumber\\
&  +\frac{1}{2\pi}%
{\displaystyle\int_{0}^{\pi}}
\left(  1-\frac{r\cos(\theta)}{\sqrt{1-r^{2}\sin^{2}\theta}}\right)
\varphi(\psi_{-}(r,\theta))d\theta\nonumber\\
&  =\frac{1}{2\pi}%
{\displaystyle\int_{0}^{\pi}}
\left(  1+\frac{r\cos(\theta)}{\sqrt{1-r^{2}\sin^{2}\theta}}\right)
\varphi(\psi_{+}(r,\theta))d\theta\nonumber\\
&  +\frac{1}{2\pi}%
{\displaystyle\int_{\pi}^{2\pi}}
\left(  1-\frac{r\cos(\theta-\pi)}{\sqrt{1-r^{2}\sin^{2}(\theta-\pi)}}\right)
\varphi(\psi_{-}(r,\theta-\pi))d\theta\nonumber\\
&  =\frac{1}{2\pi}%
{\displaystyle\int_{0}^{\pi}}
\left(  1+\frac{r\cos(\theta)}{\sqrt{1-r^{2}\sin^{2}\theta}}\right)
\varphi(\psi_{+}(r,\theta))d\theta\nonumber\\
&  +\frac{1}{2\pi}%
{\displaystyle\int_{\pi}^{2\pi}}
\left(  1+\frac{r\cos(\theta)}{\sqrt{1-r^{2}\sin^{2}\theta}}\right)
\varphi(\psi_{-}(r,\theta-\pi))d\theta\nonumber\\
&  =\frac{1}{2\pi}%
{\displaystyle\int_{0}^{2\pi}}
\left(  1+\frac{r\cos(\theta)}{\sqrt{1-r^{2}\sin^{2}\theta}}\right)
\varphi(\psi_{\ast}(r,\theta))d\theta\label{eq2002}%
\end{align}
Here we used the variable transformation $\theta\mapsto\pi+\theta$ for the
second equality and the new definition
\begin{equation}
\psi_{\ast}(r,\theta):=\left\{
\begin{array}
[c]{cl}%
\psi_{+}(r,\theta), & 0\leq\theta\leq\pi\smallskip\\
\psi_{-}(r,\theta-\pi), & \pi<\theta\leq2\pi
\end{array}
\right. \label{eq2003}%
\end{equation}
where the functions $\psi_{-}$ and $\psi_{+}$ are defined in (\ref{eq2000})
and (\ref{eq2001}). We remark, that the function $\psi_{\ast}:[0,1)\times
\lbrack0,2\pi]\mapsto\lbrack0,2\pi]$ is a continuous function which follows
from its geometric construction. We further define
\begin{equation}
k(r,\theta):=1+\frac{r\cos(\theta)}{\sqrt{1-r^{2}\sin^{2}\theta}%
}\label{eq2004}%
\end{equation}
a $2\pi$ periodic continuous function on $[0,1)\times\lbrack0,2\pi]$. If we
now go back to the general case $P=r(\cos(\alpha),\sin(\alpha))$, $\alpha
\in\lbrack0,2\pi)$, we can rotate the given boundary function $\varphi$ and
obtain
\begin{equation}
(\mathcal{E}\varphi)(P)\equiv\Phi(P)=\frac{1}{2\pi}\int_{0}^{2\pi}%
k(r,\theta)\varphi(\psi_{\ast}(r,\theta)+\alpha)d\theta\label{eq2005}%
\end{equation}
Before we study the properties of the extension operator $\mathcal{E}$, we
present two numerical examples.

To obtain $\Phi\left(  P\right)  $, we apply the trapezoidal rule to
approximate the integral in (\ref{eq2005}) or (\ref{eq25}). The number of
integration nodes should be chosen sufficiently large, although
experimentation is needed to determine an adequate choice.

\begin{example}
Consider%
\begin{equation}
\varphi\left(  \cos\theta,\sin\theta\right)  =\left(  \cos\theta-\sin
\theta+a\cos^{2}\theta,\cos\theta+\sin\theta\right)  ,\quad0\leq\theta\leq
2\pi\label{eq30}%
\end{equation}
with $0<a<1$. We choose $a=0.9$ and apply the above with $n=100$ subdivisions
for the trapezoidal rule to evaluate (\ref{eq25}). Figure \ref{fig5} shows the
mapping $\Phi$, done in the same manner as earlier with Figures \ref{fig1} and
\ref{ovals3}.
\end{example}

\begin{example}
\label{Exam_Ovals1}We consider again the ovals of Cassini region with boundary
given in (\ref{eq20}) with $a=1.5$. The mapping (\ref{eq25}) is illustrated in
Figure \ref{fig6}. However, for $a$ somewhat closer to 1, this integration
formula (\ref{eq25}) no longer produces a satisfactory $\Phi$.
\end{example}

\subsection{Properties of $\mathcal{E}\varphi$}

To study the properties of the extension operator $\mathcal{E}$, see
(\ref{eq2005}), we have to study the behavior of the functions $\psi_{\ast}$
and $k$, see (\ref{eq2003}) and (\ref{eq2004}), at the boundary $r=1$. We
start with the function $\psi_{\ast}$ and define the values of this function
for $r=1$ first:
\begin{equation}
\psi_{\ast}(1,\theta):=\left\{
\begin{array}
[c]{cl}%
0, & 0\leq\theta\leq\frac{1}{2}\pi\smallskip\\
2\theta-\pi, & \frac{1}{2}\pi\leq\theta\leq\frac{3}{2}\pi\smallskip\\
2\pi, & \frac{3}{2}\pi\leq\theta\leq2\pi
\end{array}
\right. \label{eq2006}%
\end{equation}
Because of
\[
\lim_{r\rightarrow1}1-r^{2}=0
\]
and the boundedness of the sine function, the limit
\[
\lim_{r\rightarrow1}1-r^{2}\sin^{2}(\theta)=1-\sin^{2}(\theta),
\]
is uniform for $\theta\in\left[  0,\frac{1}{2}\pi\right]  $. The uniform
continuity of the square root function implies that
\begin{align*}
\lim_{r\rightarrow1^{-}}\left(  \cos(\theta)\sqrt{1-r^{2}\sin^{2}\theta
}\right)   &  =\cos(\theta)\sqrt{1-\sin^{2}\theta}\\
&  =\cos^{2}\theta
\end{align*}
uniformly for $\theta\in\left[  0,\frac{1}{2}\pi\right]  $. Together with
similar arguments for the function $r$ we get
\[
\lim_{r\rightarrow1^{-}}\left(  r\sin^{2}\theta+\cos(\theta)\sqrt{1-r^{2}%
\sin^{2}\theta}\right)  =\sin^{2}\theta+\cos^{2}\theta=1
\]
uniformly in $\theta$. Finally we use the uniform continuity of $\arccos
(\cdot)$ to conclude that
\begin{align*}
\lim_{r\rightarrow1^{-}}\psi_{\ast}(r,\theta)  &  =\lim_{r\rightarrow1^{-}%
}\arccos\left(  r\sin^{2}\theta+\cos(\theta)\sqrt{1-r^{2}\sin^{2}\theta
}\right) \\
&  =\arccos(1)\\
&  =0\:=\psi_{\ast}(1,\theta)
\end{align*}
converges uniformly on $\left[  0,\frac{1}{2}\pi\right]  $. Because
\[
\sqrt{1-\sin^{2}(\theta)}=-\cos(\theta),\quad\quad\theta\in\left[  \tfrac
{1}{2}\pi,\pi\right]  ,
\]
we see in a similar way that
\begin{align*}
\lim_{r\rightarrow1^{-}}\psi_{\ast}(r,\theta)  &  =\lim_{r\rightarrow1^{-}%
}\arccos\left(  r\sin^{2}\theta+\cos(\theta)\sqrt{1-r^{2}\sin^{2}\theta
}\right) \\
&  =\arccos(\sin^{2}\theta+\cos(\theta)(-\cos(\theta)))\\
&  =\arccos(-\cos(2\theta))\\
&  =\arccos(\cos(2\theta-\pi))\\
&  =2\theta-\pi\\
&  =\psi_{\ast}(1,\theta)
\end{align*}
uniformly for $\theta\in\left[  \frac{1}{2}\pi,\pi\right]  $. Similar
arguments apply for $\theta\in\lbrack\pi,2\pi]$ and we finally conclude that
$\psi_{\ast}(r,\theta)$ converges uniformly to $\psi_{\ast}(1,\theta)$ as $r$
approaches $1$. This proves the next lemma.

\begin{lemma}
\label{lemma2000} The function $\psi_{\ast}$, defined by (\ref{eq2003}) and
(\ref{eq2006}), is continuous on $[0,1]\times[0,2\pi]$.
\end{lemma}

We remark that continuity on a closed interval implies uniform continuity.

Now we turn to the function $k$ defined in (\ref{eq2004}). Here we define $k$
for the value $r=1$ in the following way
\begin{equation}
k(1,\theta):=\left\{
\begin{array}
[c]{rl}%
2, & 0\leq\theta<\frac{1}{2}\pi\smallskip\\
0, & \frac{1}{2}\pi\leq\theta\leq\frac{3}{2}\pi\smallskip\\
2, & \frac{3}{2}\pi<\theta\leq2\pi
\end{array}
\right. \label{eq2007}%
\end{equation}
Obviously $k(1,\cdot)$ cannot be the uniform limit of $k(r,\cdot)$ as $r$
approaches $1$, but the following lemma holds.

\begin{lemma}
\label{lemma2001} The function $k:[0,1]\times\lbrack0,2\pi]\mapsto\lbrack
0,2]$, defined by (\ref{eq2004}) and (\ref{eq2007}), is bounded; and for every
$\delta>0,$ the function $k(r,\theta)$ approaches $k(1,\theta)$ uniformly on
$I_{\delta}$ as $r$ approaches $1$. Here
\[
I_{\delta}:=[0,2\pi]\setminus\left\{  \left(  \tfrac{1}{2}\pi-\delta,\tfrac
{1}{2}\pi+\delta\right)  \cup\left(  \tfrac{3}{2}\pi-\delta,\tfrac{3}{2}%
\pi+\delta\right)  \right\}
\]

\end{lemma}

\begin{proof}
That $k$ is bounded by $2$ follows from
\[
\sqrt{1-r^{2}\sin^{2}\theta}\geq\sqrt{1-\sin^{2}\theta}=|\cos(\theta)|
\]
and the fact that $r\in\lbrack0,1]$. The function $1/\sqrt{z}$ is uniformly
continuous on $\varepsilon\leq z\leq1$ for every $\varepsilon>0$. From the
proof of Lemma \ref{lemma2000} we know that
\[
\lim_{r\rightarrow1^{-}}1-r^{2}\sin^{2}\theta=\cos^{2}\theta
\]
uniformly for $\theta\in\lbrack0,2\pi]$. Together with the uniform continuity
of $1/\sqrt{z}$ on $[\cos^{2}(\frac{1}{2}\pi-\delta),1],$ this shows
\[
\lim_{r\rightarrow1^{-}}1+\frac{r\cos(\theta)}{\sqrt{1-r^{2}\sin^{2}\theta}%
}=1+\frac{\cos(\theta)}{|\cos(\theta)|}%
\]
uniformly on $I_{\delta}$. Remembering
\[
|\cos(\theta)|=\left\{
\begin{array}
[c]{rl}%
\cos(\theta), & \theta\in\left[  0,\frac{1}{2}\pi\right]  \cup\left[  \frac
{3}{2}\pi,2\pi\right]  \smallskip\\
-\cos(\theta), & \theta\in\left[  \frac{1}{2}\pi,\frac{3}{2}\pi\right]
\end{array}
\right.
\]
proves the Lemma.$\left.  {}\right.  $\hfill\hfill\medskip
\end{proof}

Motivated by the properties of $\psi_{\ast}$ and $k$ we now prove a more
general result for integral operators of the form (\ref{eq2005}).

\begin{lemma}
\label{lemma2002} Let $k_{1},k_{2}:[0,1]\times\lbrack0,2\pi]\mapsto\mathbb{R}
$ be bounded functions which are continuous on $[0,1)\times\lbrack0,2\pi]$.
Assume there is a finite set $E=\{\theta_{1},\ldots,\theta_{n}\}$ such that
\[
\lim_{r\rightarrow1^{-}}k_{i}(r,\theta)=k_{i}(1,\theta),\quad\quad i=1,2,
\]
uniformly on $I_{\delta}:=\{\theta\in\lbrack0,2\pi]\mid|\theta-\theta_{j}%
|\geq\delta,\,j=1,\ldots,n\}$ for every $\delta>0$. Then for a periodic
continuous function $\varphi:[0,2\pi]\mapsto\mathbb{R}$ the function
\[
\Phi(r,\alpha):=\int_{0}^{2\pi}k_{1}(r,\theta)\varphi(k_{2}(r,\theta
)+\alpha)\,d\theta
\]
is continuous on $[0,1]\times\lbrack0,2\pi]$ and $2\pi$ periodic in $\alpha$.
\end{lemma}

\begin{remark}
The above lemma will apply to each component of the function $\mathcal{E}%
\varphi$ defined in (\ref{eq2005}) with $k_{1}=k$ and $k_{2}=\psi_{\ast}$ and
$E=\{\frac{1}{2}\pi,\frac{3}{2}\pi\}$. This shows the continuity of
$\mathcal{E}\varphi$.
\end{remark}

\begin{proof}
The uniform convergence on $I_{\delta}$, $\delta>0$ arbitrary, shows that
$k_{i}(1,\cdot)$, $i=1,2$, are piecewise continuous and bounded functions on
$[0,2\pi]$, so all integrals exist. The continuity of $\Phi(r,\alpha)$ on
$[0,1)\times\lbrack0,2\pi]$ follows easily from the continuity of the
functions $k_{i}$, $i=1,2$. The periodicity follows from the periodicity of
$\varphi$ and the definition of $\Phi$. So we only need to show the continuity
of $\Phi(r,\alpha)$ on $\{1\}\times\lbrack0,2\pi]$ for example at
$(1,\overline{\alpha})$. Because of the periodicity of $\Phi(r,\cdot)$ and the
property $(\mathcal{E}\varphi)(\alpha)=(\mathcal{E}\varphi_{\alpha})(0)$,
where $\varphi_{\alpha}(\theta)=\varphi(\alpha+\theta)$ we only need to prove
the continuity for one value of $\overline{\alpha}$, for example $\alpha=\pi$.
We estimate
\begin{align*}
|\Phi(r,\alpha)-\Phi(1,\pi)|  &  =\left\vert \int_{0}^{2\pi}k_{1}%
(r,\theta)\varphi(k_{2}(r,\theta)+\alpha)-k_{1}(1,\theta)\varphi
(k_{2}(1,\theta)+\pi)\,d\theta\right\vert \\
&  \leq\left\vert \int_{0}^{2\pi}k_{1}(r,\theta)(\varphi(k_{2}(r,\theta
)+\alpha)-\varphi(k_{2}(1,\theta)+\pi))\,d\theta\right\vert \\
&  +\left\vert \int_{0}^{2\pi}(k_{1}(r,\theta)-k_{1}(1,\theta))\varphi
(k_{2}(1,\theta)+\pi))\,d\theta\right\vert
\end{align*}
Now we know that $k_{1},k_{2},$ and $\varphi$ are bounded functions, for
example
\[
|k_{1}(r,\theta)|,|k_{1}(r,\theta)|,|\varphi(\theta)|\leq M,\quad\quad M>0
\]
for all $(r,\theta)\in\lbrack0,1]\times\lbrack0,2\pi]$. So we only have to
show
\begin{equation}
\lim_{r\rightarrow1}\int_{0}^{2\pi}|\varphi(k_{2}(r,\theta)+\alpha
)-\varphi(k_{2}(1,\theta)+\pi)|\,d\theta=0\label{eq2008}%
\end{equation}
and%
\begin{equation}
\lim_{r\rightarrow1}\int_{0}^{2\pi}|k_{1}(r,\theta)-k_{1}(1,\theta
)|\,d\theta=0\label{eq2009}%
\end{equation}
We start with the first limit. Given an $\varepsilon>0$ we choose $\delta>0 $
small enough such that
\begin{equation}
\int_{\lbrack0,2\pi]\setminus I_{\delta}}2M\;d\theta=\frac{\varepsilon}%
{2}\label{eq2010}%
\end{equation}
Now we observe that $\varphi$ is uniformly continuous on $\mathbb{R}$ because
it is continuous and periodic. So there is a $\omega>0$ such that
\begin{equation}
|\varphi(x)-\varphi(y)|\leq\frac{\varepsilon}{4\pi}\label{eq2011}%
\end{equation}
if $|x-y|\leq\omega$. We also know that $k_{2}(r,\cdot)$ converges uniformly
on $I_{\delta}$ to $k_{2}(1,\cdot)$, so there is a $r_{0}\in(0,1)$, such that
\[
|k_{2}(r,\theta)-k_{2}(1,\theta)|\leq\frac{\omega}{2}%
\]
for all $r\geq r_{0}$ and $\theta\in I_{\delta}$. If furthermore $|\alpha
-\pi|\leq\omega/2$, we conclude that
\begin{align*}
|(k_{2}(r,\theta)+\alpha)-(k_{2}(1,\theta)+\pi)|  &  \leq|k_{2}(r,\theta
)-k_{1}(1,\theta)|+|\alpha-\pi|\\
&  \leq\omega
\end{align*}
which by (\ref{eq2011}) implies
\begin{equation}
|\varphi(k_{2}(r,\theta)+\alpha)-\varphi(k_{2}(1,\theta)+\pi)|\leq
\frac{\varepsilon}{4\pi}\label{eq2012}%
\end{equation}
for all $(r,\alpha)\in\lbrack r_{0},1]\times\lbrack\pi-\omega/2,\pi+\omega/2]
$ and $\theta\in I_{\delta}$. Combining (\ref{eq2010}) and (\ref{eq2012}) we
can estimate
\begin{align*}
&  \int_{0}^{2\pi}|\varphi(k_{2}(r,\theta)+\alpha)-\varphi(k_{2}(1,\theta
)+\pi)|\,d\theta\\
&  \leq\int_{[0,2\pi]\setminus I_{\delta}}|\varphi(k_{2}(r,\theta
)+\alpha)-\varphi(k_{2}(1,\theta)+\pi)|\,d\theta\\
&  +\int_{I_{\delta}}|\varphi(k_{2}(r,\theta)+\alpha)-\varphi(k_{2}%
(1,\theta)+\pi)|\,d\theta\\
&  \leq\int_{[0,2\pi]\setminus I_{\delta}}2M\;d\theta+\int_{I_{\delta}}%
\frac{\varepsilon}{4\pi}\,d\theta\\
&  \leq\frac{\varepsilon}{2}+2\pi\cdot\frac{\varepsilon}{4\pi}=\varepsilon
\end{align*}
for all $(r,\alpha)\in\lbrack r_{0},1]\times\lbrack\pi-\omega/2,\pi+\omega
/2]$,which proves (\ref{eq2008}).

To prove (\ref{eq2009}) we again choose an arbitrary $\varepsilon>0$ and
choose $\delta>0$ such that (\ref{eq2010}) is true. Now the uniform
convergence of $k_{1}(r,\cdot)$ to $k_{1}(1,\cdot)$ on $I_{\delta}$ proves the
existence of a $r_{1}\in(0,1)$ such that
\begin{equation}
|k_{1}(r,\theta)-k_{1}(1,\theta)|\leq\frac{\varepsilon}{4\pi}\label{eq2013}%
\end{equation}
for all $(r,\theta)\in\lbrack r_{1},1]\times I_{\delta}$. Using (\ref{eq2010})
and (\ref{eq2013}) we estimate
\begin{align*}
\int_{0}^{2\pi}|k_{1}(r,\theta)-k_{1}(1,\theta)|\;d\theta &  =\int%
_{[0,2\pi]\setminus I_{\delta}}|k_{1}(r,\theta)-k_{1}(1,\theta)|\,d\theta\\
&  +\int_{I_{\delta}}|k_{1}(r,\theta)-k_{1}(1,\theta)|\,d\theta\\
&  \leq\int_{[0,2\pi]\setminus I_{\delta}}2M\,d\theta+\int_{I_{\delta}}%
\frac{\varepsilon}{4\pi}\,d\theta\\
&  \leq\frac{\varepsilon}{2}+2\pi\cdot\frac{\varepsilon}{4\pi}=\varepsilon
\end{align*}
for all $r\in\lbrack r_{1},1]$. This proves (\ref{eq2009}).$\left.  {}\right.
$\hfill\hfill\medskip
\end{proof}

Now we state the results about the extension operator $\mathcal{E}$.

\begin{theorem}
\label{theorem2000} Let $\varphi:\partial B_{2}\mapsto\mathbb{R}^{2}$, be a
continuous function, then $\Phi(P)=(\mathcal{E}\varphi)(r,\alpha)$,
$P\in\overline{B_{2}}$, see (\ref{eq2005}), is continuous function on
$\overline{B_{2}}$ and
\begin{align}
\Phi|_{\partial B_{2}}  &  = \varphi\label{eq2014}%
\end{align}

\end{theorem}

\begin{proof}
In Lemma \ref{lemma2000} and Lemma \ref{lemma2001} we have shown that the
functions $k$ and $\psi_{\ast}$ in (\ref{eq2005}) satisfy the assumptions of
Lemma \ref{lemma2002}. So the continuity of $\Phi(P)$ follows from Lemma
\ref{lemma2002}. For $P\in\partial B_{2}$ the polar coordinates of $P$ are
given by $(r,\alpha)=(1,\alpha)$, $\alpha\in\lbrack0,2\pi]$, so we get with
(\ref{eq2006}) and (\ref{eq2007})
\begin{align*}
\Phi(P)  &  =\frac{1}{2\pi}\int_{0}^{2\pi}k(1,\theta)\varphi(\psi_{\ast
}(1,\theta)+\alpha)\,d\theta\\
&  =\frac{1}{2\pi}\left(  \int_{0}^{\pi/2}2\varphi(0+\alpha)\;d\theta
+\int_{3\pi/2}^{2\pi}2\varphi(2\pi+\alpha)\,d\theta\right) \\
&  =\frac{1}{2\pi}(\pi\varphi(\alpha)+\pi\varphi(2\pi+\alpha))\\
&  =\varphi(\alpha)=\varphi(P)
\end{align*}
because of the $2\pi$ periodicity of $\varphi$.$\left.  {}\right.  $%
\hfill\hfill\medskip
\end{proof}

\begin{corollary}
\label{corollary2000} Let $\Omega\subset\mathbb{R}^{2}$ be a domain with
boundary $\partial\Omega$ and $\varphi:\partial B_{2}\mapsto\partial\Omega$ be
a continuous parametrization of the boundary. Then the function $\mathcal{E}%
\varphi$, defined in (\ref{eq2005}), maps $\overline{B_{2}}$ onto
$\overline{\Omega}$.
\end{corollary}

\begin{proof}
Theorem \ref{theorem2000} implies that $\mathcal{E}\varphi:\overline{B_{2}%
}\mapsto\mathbb{R}^{2}$ is continuous and that $(\mathcal{E}\varphi)(\partial
B_{2})=\partial\Omega$. We assume that the parametrization $\varphi$ moves
along the boundary of $\Omega$ in the positive direction. For $y\in\Omega$ we
then have
\[
\deg(\mathcal{E}\varphi,y)=1
\]
where $\deg$ is the mapping degree; see \cite[Chapter 12]{ZEI}. But this
implies that there is at least one $x\in B_{2}$ such that $(\mathcal{E}%
\varphi)(x)=y$.
\end{proof}

\begin{theorem}
\label{theorem2001} Let $\Omega\subset\mathbb{R}^{2}$ be a convex domain with
boundary $\partial\Omega$ and $\varphi:\partial B_{2}\mapsto\partial\Omega$ be
a continuous parametrization of the boundary. Then $(\mathcal{E}\varphi
)(B_{2})\subset\overline{\Omega}$.
\end{theorem}

\begin{proof}
We have to show $(\mathcal{E}\varphi)(P)\in\overline{\Omega}$ for every $P\in
B_{2}$. We use the first equation in formula (\ref{eq2002})
\begin{align*}
\Phi(P)  &  =\frac{1}{2\pi}\int_{0}^{\pi}\left(  1+\frac{r\cos(\theta)}%
{\sqrt{1-r^{2}\sin^{2}(\theta)}}\right)  \varphi(\psi_{+}(r,\theta))\\
&  +\left(  1-\frac{r\cos(\theta)}{\sqrt{1-r^{2}\sin^{2}(\theta)}}\right)
\varphi(\psi_{-}(r,\theta))\,d\theta\\
&  =\lim_{N\rightarrow\infty}\frac{1}{N}\sum_{j=0}^{N}\left(  \frac{1}%
{2}+\frac{r\cos(\theta_{j})}{2\sqrt{1-r^{2}\sin^{2}(\theta_{j})}}\right)
\varphi(\psi_{+}(r,\theta_{j}))\\
&  +\left(  \frac{1}{2}-\frac{r\cos(\theta_{j})}{2\sqrt{1-r^{2}\sin^{2}%
(\theta_{j})}}\right)  \varphi(\psi_{-}(r,\theta_{j}))
\end{align*}
where $\theta_{j}=\pi j/N$ and we further assumed again that $P$ is on the
positive real axis to simplify the notation. Here we have used the fact that
the integral is the limit of Riemann sums. Each term of the sum is a convex
combination of two elements of $\overline{\Omega}$ and therefore in
$\overline{\Omega}$. But the sum itself is a convex combination, so the sum is
an element of $\overline{\Omega}$. Finally $\overline{\Omega}$ is closed, so
$\Phi(P)\in\overline{\Omega}$.$\left.  {}\right.  $\hfill\hfill\medskip
\end{proof}

The two last results imply that for a convex domain $\Omega$ we get
$(\mathcal{E}\varphi)(\overline{B_{2}})=\overline{\Omega}$, but there is still
the possibility that $\mathcal{E}(\varphi)$ is not injective. Our numerical
examples seem to indicate that the function is injective for convex $\Omega$,
but we have no proof. For non-convex regions, it works for some but not
others. It is another option in a toolkit of methods for producing the mapping
$\Phi$.

The integration-based formula (\ref{eq25}) can be extended to three
dimensions. Given%
\[
\varphi:\partial B_{3}\underset{onto}{\overset{1-1}{\longrightarrow}}%
\partial\Omega,
\]
define the interpolation formula $\varphi_{\ast}\left(  \theta,\omega\right)
$ as before in (\ref{eq24}), with $\left(  \theta,\omega\right)  $ the
spherical coordinates of a direction vector $\boldsymbol{\eta}$ through a
given point $P\in B_{3}$. Then define
\begin{equation}
\Phi\left(  P\right)  =\frac{1}{2\pi}\int_{0}^{2\pi}\int_{0}^{\pi/2}%
\varphi_{\ast}\left(  \theta,\omega\right)  \sin\omega\,d\omega\,d\theta
\label{eq2020}%
\end{equation}
A proof of the generalization of Corollary \ref{corollary2000} can be given
along the same line as given above for the disk $B_{2}$.
%

\begin{figure}[tb]%
\centering
\includegraphics[
height=3in,
width=3.9998in
]%
{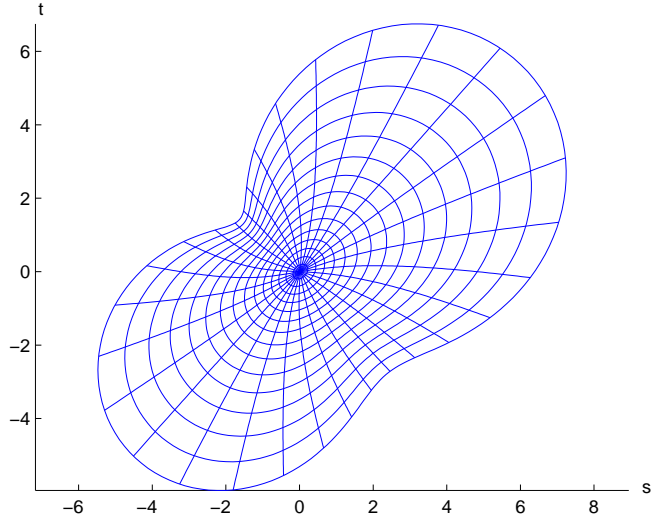}%
\caption{The initial mapping $\protect\widetilde{\Phi}$ for Example
\ref{Exam_IT1} with $a=5$, \ based on (\ref{e40})}%
\label{fig7}%
\end{figure}
%

\begin{figure}[tb]%
\centering
\includegraphics[
height=3in,
width=3.9998in
]%
{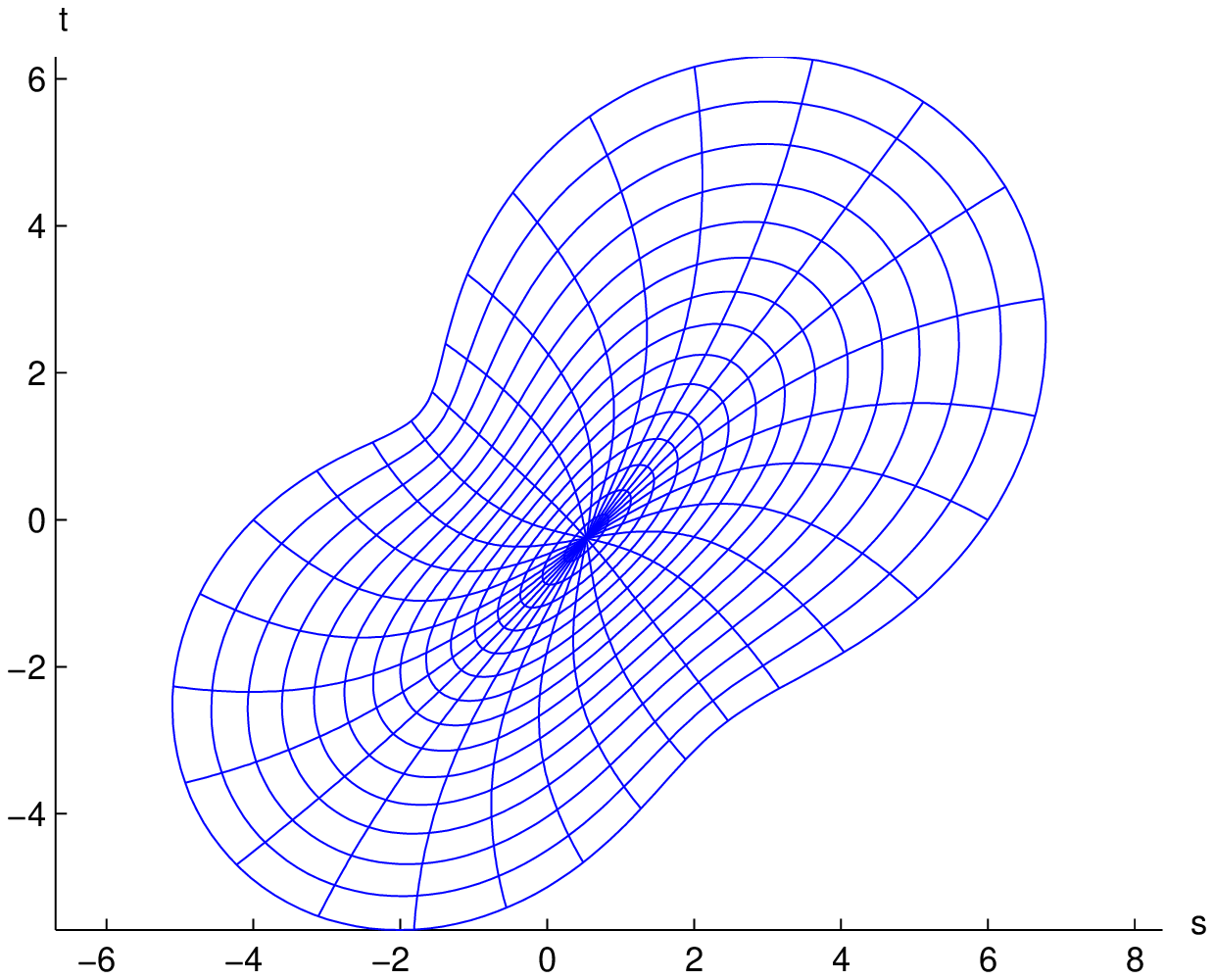}%
\caption{The mapping $\Phi$ for Example \ref{Exam_IT1} with $a=5$, obtained
using iteration}%
\label{fig8}%
\end{figure}

\section{Iteration methods\label{iteration}}

Some of the methods discussed in \S 2 lead to a mapping $\Phi$ in which
$\det\left(  D\Phi\left(  x\right)  \right)  $ has a large variation as $x$
ranges over the unit ball $B_{d}$, especially those methods based on using
using the $C^{\infty}$-function $T\left(  r,\kappa\right)  $ of (\ref{e40}).
\ We seek a way to improve on such a mapping, to obtain a mapping in which
$\det\left(  D\Phi\left(  x\right)  \right)  $ has a smaller variation over
$B_{d} $. We continue to look at only the planar problem, while keeping in
mind the need for a method that generalizes to higher dimensions. In this
section we introduce an iteration method to produce a mapping $\Phi$ with each
component a multivariate polynomial over $B_{2}$.

Assume we have an initial guess for our mapping, in the form of a polynomial
of degree $n$,%
\begin{equation}
\Phi_{n}^{(0)}\left(  x\right)  \equiv\sum_{j=1}^{N_{n}}\alpha_{j}^{(0)}%
\psi_{j}(x),\quad\quad x\in B_{2}\label{F1}%
\end{equation}
We want to calculate an `improved' value for $\Phi_{n}^{(0)}$, call it
$\Phi_{n}$.

The coefficients $\alpha_{j}^{(0)}\in\mathbb{R}^{2}$. The functions $\left\{
\psi_{1},\dots,\psi_{N_{n}}\right\}  $ are chosen to be a basis for $\Pi_{n}$,
the polynomials of degree $\leq n$. and we require them to be orthonormal with
respect to the inner product $\left(  \cdot,\cdot\right)  $ associated with
$L^{2}\left(  B_{2}\right)  $. Note that $N_{n}=\dim\left(  \Pi_{n}\right)
=\frac{1}{2}\left(  n+1\right)  \left(  n+2\right)  $. As basis functions
$\left\{  \psi_{j}\right\}  $ in our numerical examples, we use the `ridge
polynomials' of Logan and Shepp \cite{Loga}, an easy basis to define and
calculate; also see \cite[\S 4.3.1]{atkinson-han2011}.

We use an iterative procedure to seek an approximation%
\begin{equation}
\Phi_{n}\left(  x\right)  =\sum_{j=1}^{N_{n}}\alpha_{n,j}\psi_{j}(x)\label{F3}%
\end{equation}
of degree $n$ that is an improvement in some sense on $\Phi_{n}^{(0)}$. \ The
degree $n$ used in defining $\Phi_{n}^{\left(  0\right)  }$, and also in
defining our improved value $\Phi_{n}$, will need to be sufficiently large;
and usually, $n$ must be determined experimentally.

The coefficients $\left\{  \alpha_{j}^{(0)}\right\}  $ are normally generated
by numerical integration of the Fourier coefficients $\left\{  \alpha
_{j}^{(0)}\right\}  $,
\begin{equation}
\alpha_{j}^{(0)}=\left(  \widetilde{\Phi},\psi_{j}\right)  ,\label{F2}%
\end{equation}
where $\widetilde{\Phi}$ is generated by one of the methods discussed in
\S \S \ref{constructions},\ref{integinterp}. The quadrature used is
\begin{equation}
\int_{B_{2}}g(x,y)\,dx\,dy\approx\frac{2\pi}{2p+1}\sum_{l=0}^{p}\sum
_{m=0}^{2p}\omega_{l}r_{l}\widehat{g}\left(  r_{l},\frac{2\pi\,m}{2p+1}\right)
\label{e106}%
\end{equation}
where $\widehat{g}\left(  r,\theta\right)  \equiv g\left(  r\cos\theta
,r\sin\theta\right)  $. Here the numbers $\left\{  \omega_{l}\right\}  $ are
the weights of the $\left(  p+1\right)  $-point Gauss-Legendre quadrature
formula on $[0,1]$, and the nodes $\left\{  r_{l}\right\}  $ are the
corresponding zeros the the degree $p+1$ Legendre polynomial on $\left[
0,1\right]  $. \ This formula is exact if $g$ is a polynomial of degree
$\leq2p+1 $; see \cite[\S 2.6]{stroud}.\ 

We need to require that our mapping will agree with $\varphi$ on $S^{1}$, at
least approximately. \ To this end, choose a formula $q_{n}$ for the number of
points on $S^{1}$ at which to match $\Phi_{n}$ with the function $\varphi$ and
then choose $\left\{  z_{1},\dots,z_{q_{n}}\right\}  $ on $S^{1} $. Require
$\Phi_{n}$ to satisfy%
\begin{equation}
\Phi_{n}\left(  z_{j}\right)  =\varphi\left(  z_{j}\right)  ,\quad\quad
j=1,\dots,q_{n}\label{F4a}%
\end{equation}
which imposes implicitly $q_{n}$ conditions on the coefficients $\left\{
\alpha_{n,j}\right\}  $. If $\varphi$ is a trigonometric polynomial of degree
$m$, and if $n\geq m$ with $q_{n}=2n+1$, then (\ref{F4a}) will imply that
$\left.  \Phi_{n}\right\vert _{S^{1}}=\varphi$ over $\partial\Omega$. Our
numerical examples all use this latter choice of $q_{n}$.

Next, choose a function $\mathcal{F}\left(  \alpha\right)  $, $\alpha=\left[
\alpha_{1},\dots,\alpha_{N}\right]  ^{\text{T}}$, and seek to calculate
$\alpha$ so as to minimize $\mathcal{F}(\alpha)$ subject to the above
constraints (\ref{F4a}). How should $\mathcal{F}$ be chosen? To date, the most
successful choice experimentally has been $\mathcal{F}\left(  \alpha\right)
=\Lambda\left(  \Phi_{n}\right)  ,$ defined earlier in (\ref{eq10}).
%

\begin{figure}[tb]%
\centering
\includegraphics[
height=3in,
width=3.9998in
]%
{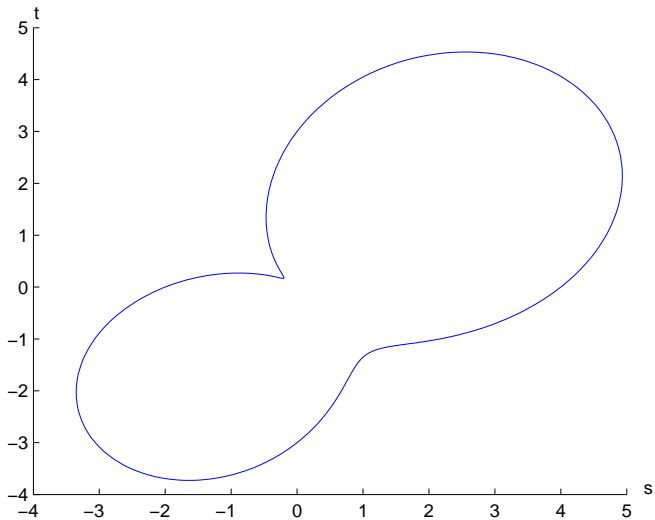}%
\caption{The boundary for the starlike region with $\rho=3+\cos\theta
+2\sin2\theta$}%
\label{fig13}%
\end{figure}
%

\begin{figure}[tb]%
\centering
\includegraphics[
height=3in,
width=3.9998in
]%
{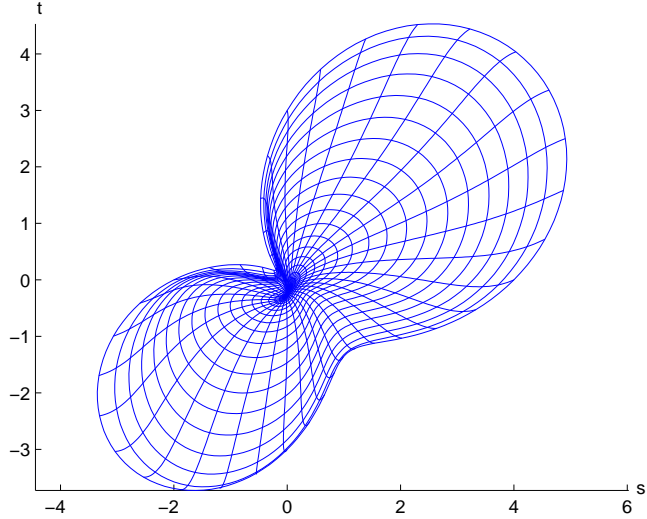}%
\caption{The boundary mapping $\Phi$ for the starlike region \newline with
$\rho=3+\cos\theta+2\sin2\theta$}%
\label{fig14}%
\end{figure}
%

\begin{figure}[tb]%
\centering
\includegraphics[
height=3in,
width=3.9998in
]%
{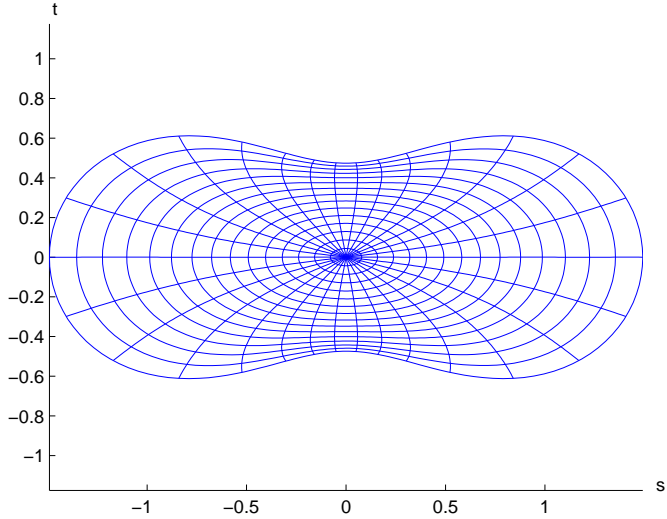}%
\caption{The boundary mapping $\Phi$ for the starlike region with $\rho$ from
(\ref{eq20}) with $a=1.5$}%
\label{fig16}%
\end{figure}
\begin{figure}[tb]%
\centering
\includegraphics[
height=3in,
width=3.9998in
]%
{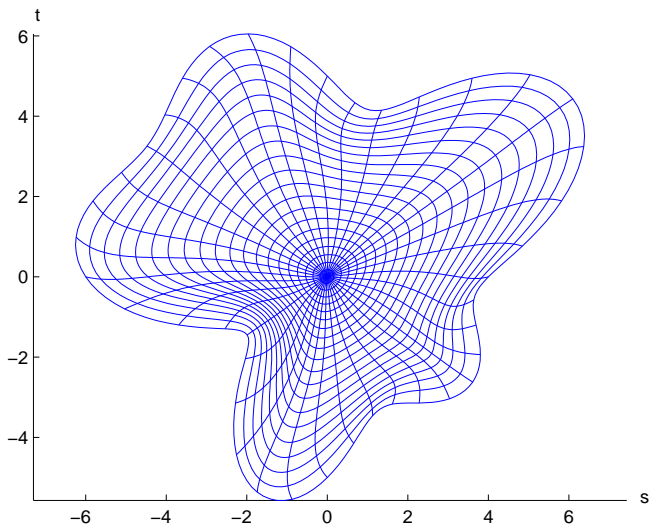}%
\caption{The optimal mapping $\Phi_{7}$ for the starlike region with
$\protect\widehat{\rho}\left(  \theta\right)  =5+\sin\theta+\sin3\theta
-\cos5\theta$}%
\label{amoeba_map}%
\end{figure}

\subsection{The iteration algorithm\label{sec4.1}}

Using the constraints (\ref{F4a}) leads to the system%
\begin{equation}
A\alpha=\mathbf{\varphi},\label{F4}%
\end{equation}%
\[
A=\left[
\begin{array}
[c]{ccc}%
\psi_{1}(z_{1}) & \cdots & \psi_{N}(z_{1})\\
\vdots &  & \vdots\\
\psi_{1}(z_{q}) & \cdots & \psi_{N}(z_{q})
\end{array}
\right]  ,\quad\quad\mathbf{\varphi=}\left[
\begin{array}
[c]{c}%
\varphi_{1}(z_{1})\\
\vdots\\
\varphi_{1}(z_{q})
\end{array}
\right]
\]
Because $\left.  \Phi_{n}\right\vert _{S^{1}}$ is a trigonometric polynomial
of degree $n$, it is a bad idea to have $q_{n}>2n+1$. The maximum row rank of
$A$ can be at most $2n+1$. Let $\left\{  z_{1},\dots,z_{q}\right\}  $ denote
$q_{n}$ evenly spaced points on $S^{1}$. We want to minimize $\mathcal{F}%
\left(  \alpha\right)  $ subject to the constraints (\ref{F4}).

We turn our constrained minimization problem into an unconstrained problem.
Let $A=UCV$ be the singular value decomposition of $A$; $U$ is an orthogonal
matrix of order $q$, $V$ is an orthogonal matrix of order $N$, and $C$ is a
`diagonal matrix' of order $q\times N$. The constraints (\ref{F4}) can be
written as
\begin{equation}
CV\alpha=U^{\text{T}}\mathbf{\varphi}\label{F6}%
\end{equation}
Introduce a new variable $\beta=V\alpha$, or $\alpha=V^{\text{T}}\beta$. Then
$C\beta=U^{\text{T}}\mathbf{\varphi}$ and we can solve explicitly for
$\gamma=[\beta_{1},\dots,\beta_{q}]^{\text{T}}$. Implicitly this assumes that
$A$ has full rank. Let $\delta=[\beta_{q+1},\dots,\beta_{N}]^{\text{T}}$,
$\beta=\left[  \gamma^{\text{T}},\delta^{\text{T}}\right]  ^{\text{T}} $.
\ \ Then introduce
\begin{equation}
G\left(  \delta\right)  =\mathcal{F}\left(  \alpha\right) \label{F7}%
\end{equation}
using $\alpha=V^{\text{T}}\beta$ and the known values of $\gamma$. We use our
initial $\left\{  \alpha_{j}^{(0)}\right\}  $ in (\ref{F1}) to generate the
initial value for $\beta$ and thus for $\delta$.

The drawback to this iteration method is the needed storage for the $q\times
N$ matrix $A$ and the matrices produced in its singular value decomposition.
In the following numerical examples, we minimize $G$ using the \textsc{Matlab}
program \emph{fminunc} for unconstrained minimization problems.

\begin{example}
\label{Exam_IT1}Recall Example \ref{Exam1} with $a=5$. Generate an initial
mapping $\widetilde{\Phi}$ using (\ref{e40}) with $\kappa=.5$, $\omega=1.0 $.
Next, generate an initial polynomial (\ref{F1}) of degree $n=3$, using
numerical integration of the Fourier coefficients $\left\{  \alpha_{j}%
^{(0)}\right\}  $ of (\ref{F2}). We then use the above iteration method to
obtain an improved mapping. Figure \ref{fig7} shows the initial mapping
$\widetilde{\Phi}$, and Figure \ref{fig8} shows the final mapping $\Phi_{n}$
obtained by the iteration method. With the final mapping, we have $\left.
\Phi_{n}\right\vert _{S^{1}}=\varphi$ to double precision rounding accuracy,
and%
\[
\Lambda\left(  \Phi\right)  =6.21
\]
Compare the latter to $\Lambda\left(  \Phi\right)  =100.7$ for the mapping in
Example \ref{Exam1}. \ 
\end{example}

\begin{example}
\label{Exam_IT2}Consider again the starlike region using (\ref{e20}) of
Example \ref{Exam1}, but now with $a=3$. The harmonic mapping of
\S \ref{HarmMap} failed in this case to produce a 1-1 mapping. In fact, the
boundary is quite ill-behaved in the neighborhood of $\left(  -0.2,0.2\right)
$, being almost a corner; see Figure \ref{fig13}. In this case we needed $n=7
$, with this smallest sufficient degree determined experimentally. To generate
the initial guess $\widetilde{\Phi}$, we used (\ref{e40}) with $\left(
\kappa,\omega\right)  =\left(  0.5,0.1\right)  $. For the initial guess,
$\Lambda\left(  \Phi_{7}^{(0)}\right)  \doteq840.$ We iterated first with the
\textsc{Matlab} program \emph{fminunc}. When it appeared to converge, we used
the resulting minimizer as an initial guess with a call to the \textsc{Matlab}
routine \emph{fminsearch}, which is a Nelder-Mead search method. When it
converged, its minimizer was used again as an initial guess, returning to a
call on \emph{fminunc}. Figure \ref{fig14} shows the final mapping $\Phi_{7}$
obtained with this repeated iteration. For the Jacobian matrix, $\Lambda
\left(  \Phi_{7}\right)  \doteq177.9$, further illustrating the ill-behaviour
associated with this boundary. As before, $\left.  \Phi\right\vert _{S^{1}%
}=\varphi$ to double precision rounding accuracy.
\end{example}

\begin{example}
\label{Exam_IT4}Consider again the ovals of Cassini region with boundary given
in (\ref{eq20}) with $a=1.5$. As our initial mapping $\widetilde{\Phi}$, we
use the interpolating integration-based mapping of (\ref{eq25}), illustrated
in Figure \ref{fig6}. We produce the initial guess for the coefficients
$\left\{  \alpha_{j}^{(0)}\right\}  $ of (\ref{F2}) by using numerical
integration. Unlike the preceding three examples, the boundary mapping
$\varphi$ is not a trigonometric polynomial, and thus the interpolating
conditions of (\ref{F4a}) will not force $\left.  \Phi_{n}\right\vert _{S^{1}%
}$ to equal $\varphi$ over $\partial\Omega$. For that reason, we use a higher
degree than with the preceding examples, choosing $n=16$. Figure \ref{fig16}
shows the resulting mapping $\Phi$. With this mapping, $\Lambda\left(
\Phi\right)  =26.11$. On the boundary,%
\[
\max_{x\in S^{1}}\left\vert \Phi\left(  x\right)  -\varphi\left(  x\right)
\right\vert =2.61E-4
\]
showing the mapping does not depart far from the region $\overline{\Omega}$.
\end{example}

\begin{example}
Consider the starlike domain with
\[
\widehat{\rho}\left(  \theta\right)  =5+\sin\theta+\sin3\theta-\cos
5\theta,\quad\quad0\leq\theta\leq2\pi.
\]
in (\ref{eq6})-(\ref{eq7}). Using the degree $n=7$ and the inital mapping
$\widetilde{\Phi}$ based on (\ref{e40}) with $\left(  \kappa,\omega\right)
=\left(  0.2,\,1.4\right)  $, we obtained the mapping illustrated in Figure
\ref{amoeba_map}. The minimum value obtained was $\Lambda\left(  \Phi
_{7}\right)  \doteq6.63$. As a side-note of interest, the iteration converged
to a value of $\Lambda\left(  \Phi\right)  $ that varied with the initial
choice of $\left(  \kappa,\omega\right)  $. We have no explanation for this,
other than to say that the objective function $\Lambda\left(  \Phi\right)  $
appears to be ill-behaved in some sense that we do not yet understand.
\end{example}

\subsection{An energy method}

In this section we present a second iteration method, one based on a different
objective function. Instead of $\Lambda$, see (\ref{eq10}), we use
\begin{align}
\widetilde{\Lambda}(\Phi_{n})  & \equiv\sum_{i=1}^{K_{1}}\sum_{\substack{j=1
\\i\neq j}}^{K_{1}}\frac{1}{\Vert\Phi_{n}(\xi_{i})-\Phi_{n}(\xi_{j})\Vert
_{2}^{\alpha}}\nonumber\\
& +\sum_{i=1}^{K_{1}}\sum_{j=1}^{L_{1}}\frac{1}{\Vert\Phi_{n}(\xi_{i}%
)-\Phi_{n}(\zeta_{j})\Vert_{2}^{\alpha}}.\label{eq4200}%
\end{align}
We again impose the interpolation conditions given in (\ref{F4a}); and the
free parameters are given by $\delta$, see (\ref{F7}). First we explain the
definition of the points $\xi_{i}$ and $\zeta_{j}$ appearing in formula
(\ref{eq4200}). The points $\xi_{i}$ are located inside the unit disk and are
elements of a rectangular grid
\[
\{\xi_{i}\mid i=1,\ldots,K_{1}\}=\left(  \frac{1}{k_{1}}\mathbb{Z}^{2}\right)
\cap B_{2};
\]
the density of the grid is determined by $k_{1}>0$. The points $\zeta_{j}$ are
located on the unit circle and distributed uniformly
\[
\{\zeta_{j}\mid j=1,\ldots,L_{1}\}=\left\{  \left(  \cos\left(  \frac{2\pi
j}{L_{1}}\right)  ,\sin\left(  \frac{2\pi j}{L_{1}}\right)  \right)  \mid
j=1,\ldots,L_{1}\right\}
\]
$L_{1}\in\mathbb{N}$. Furthermore the function $\widetilde{\Lambda}$ contains
the parameter $\alpha>0$. So in addition to the dimension $n$ of the trial
space for $\Phi_{n}$, this method uses four parameters: $q_{n}$, the number of
interpolation points along the boundary; $k_{1}$, which determines the grid
density inside the unit disk; $L_{1}$, the number of points along the
boundary; and $\alpha,$ the exponent in formula (\ref{eq4200}).

The motivation for the function $\widetilde{\Lambda}$ is the following. We
start with an equally distributed set of points in the unit disk, $\{\xi
_{i}\mid i=1,\ldots,K_{1}\}$ and we try to force the mapping $\Phi_{n}$ to
distribute these points as uniformly as possible in the new domain $\Omega$.
One can think of charged particles which repel each other with a certain
force. If this force is generated by the potential $r^{-\alpha}$ then the
first term in formula (\ref{eq4200}) is proportional to the energy of the
charge distribution $\{\Phi_{n}(\xi_{i})\mid i=1,\ldots,K_{1}\}$. When we go
back to our original goal of creating a mapping $\Phi$ which is injective, we
see that this is included in this fuctional because the energy becomes
infinite if two particles are moved closer.

The second goal for our mapping is that $\Phi_{n}(B_{2})\subset\Omega$, to
enforce this condition we use a particle distribution along the boundary of
$\Omega$ given by $\{\Phi_{n}(\zeta_{j})\mid j=1,\ldots,L_{1}\}$. These
charges will repel the charges $\{\Phi_{n}(\xi_{i})\mid i=1,\ldots,K_{1}\}$
away from the boundary. The energy associated with the interaction between the
interior points and the boundary points gives us the second term in formula
(\ref{eq4200}).

So we can consider the algorithm to minimize the function $\widetilde{\Lambda
}$ as an attempt to minimize the energy of a particle distribution in $\Omega
$. This should also guarantee that the mapping $\Phi_{n}$ has a small value
for the function $\Lambda$, because the original points $\{\xi_{i}\mid
i=1,\ldots,K_{1}\}$ are equally distributed.

In our numerical experiments we used $\alpha=2$, so the function
$\widetilde{\Lambda}(\Phi_{n})$ is differentiable as a function of the
parameter $\delta$. Furthermore we adjust $k_{1}\in\mathbb{N}$ in such a way
that $K_{1}\approx N_{n}$ and we choose $L_{1}\sim k_{1}$. For the parameter
$q_{n}$ we use the same value as in \S \ref{sec4.1}.

\begin{example}
Consider the starlike domain defined in (\ref{e20}) with $a=5$ again. We use
$n=3$, $\alpha=2$, $K_{1}=177$, $L_{1}=160$. To minimize the function
$\widetilde{\Lambda}$ we use the BFGS method, see \cite{Nocedal}. Figures
\ref{energy1} and \ref{energy1a} show a rectangular grid in the unit disc and
its image under the mapping $\Phi_{3}^{(0)}$. For the initial guess we have
$\widetilde{\Lambda}(\Phi_{3}^{(0)})\approx11500$ and $\Lambda(\Phi_{3}%
^{(0)})\approx29$. For the final mapping $\Phi_{3}$ we obtain
$\widetilde{\Lambda}(\Phi_{3})\approx7930$ and $\Lambda(\Phi_{3})\approx10$.
This shows that the function $\widetilde{\Lambda}$ implicitly also minimizes
the function $\Lambda$. Figure \ref{energy2} shows the image of the final
mapping $\Phi$.
\end{example}

\begin{center}%
\begin{figure}[tb]%
\centering
\includegraphics[
height=3in,
width=3.9998in
]%
{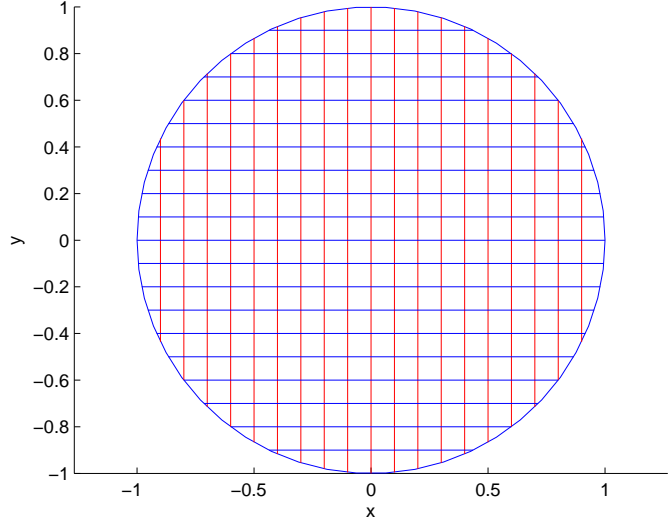}%
\caption{A grid on the unit disk}%
\label{energy1}%
\end{figure}
%

\begin{figure}[tb]%
\centering
\includegraphics[
height=3in,
width=3.9998in
]%
{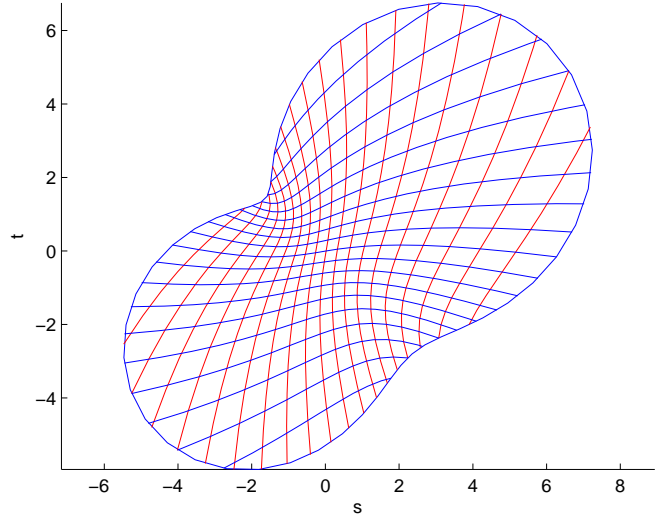}%
\caption{The image of the grid in Figure \ref{energy1} under the mapping
$\Phi_{3}^{(0)}$ for the domain given in (\ref{e20}).}%
\label{energy1a}%
\end{figure}
%

\begin{figure}[tb]%
\centering
\includegraphics[
height=3in,
width=3.9998in
]%
{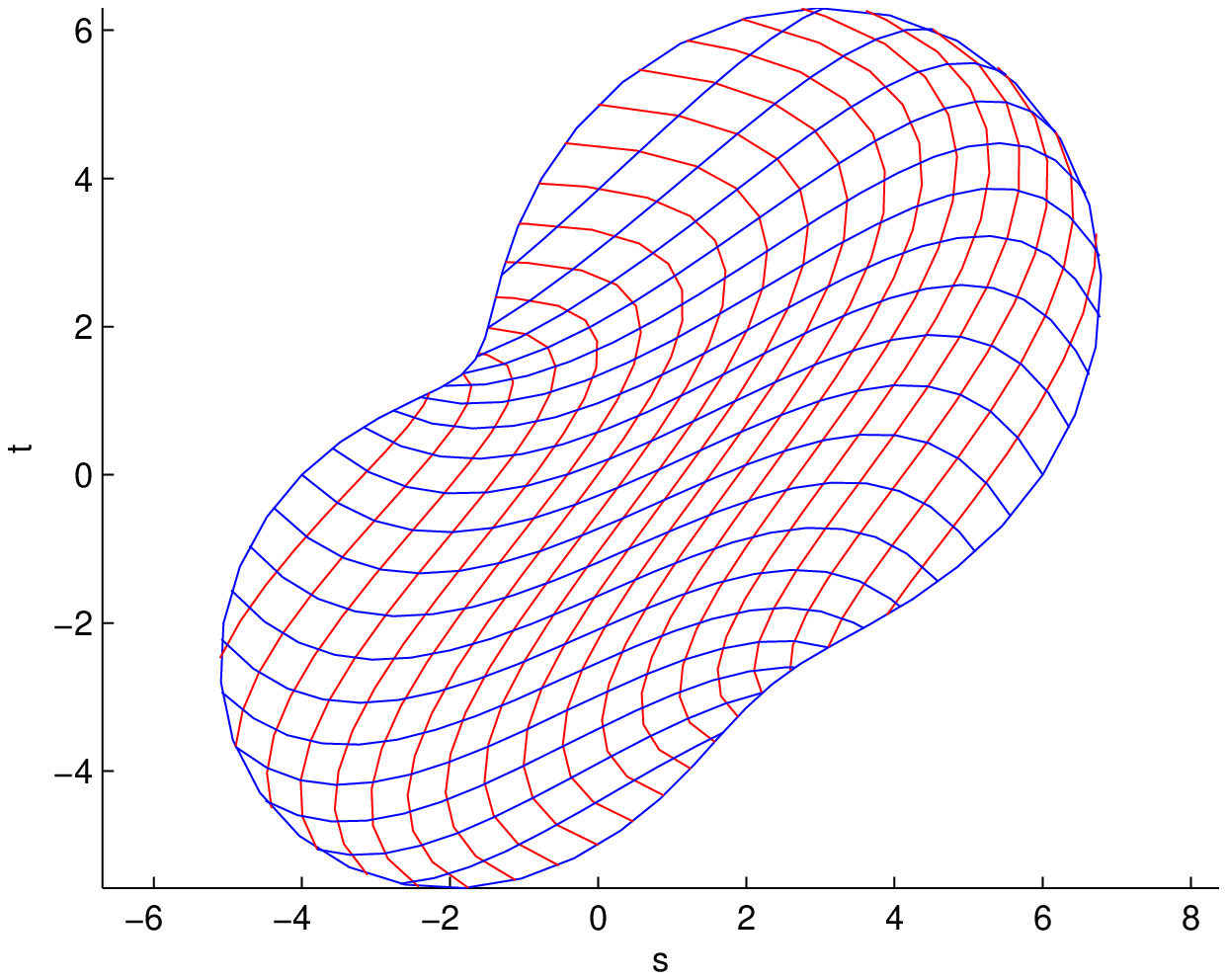}%
\caption{The image of the grid in Figure \ref{energy1} under the final mapping
$\Phi_{3}$.}%
\label{energy2}%
\end{figure}

\end{center}

\section{Mapping in three dimensions\label{3D}}

In this section we describe an algorithm to construct an extension $\Phi
_{n}:B_{3}\mapsto\Omega$ for a given function $\varphi:S^{2}\mapsto
\partial\Omega$. We assume that $\Omega$ is starlike with respect to the
origin. The three dimensional case differs from the algorithm described in
\S \ref{iteration} in several ways. The dimension of $\Pi_{n}$ of the
polynomials of maximal degree $n$ is given by $N_{n}={\binom{n}{3}}$, so any
optimization algorithm has to deal with a larger number of degrees of freedom
for a given $n$ when compared to the two dimensional case. Whereas in the two
dimensional case a plot of $\boldsymbol{\Phi}_{n}(B_{2})$ reveals any problems
of the constructed $\boldsymbol{\Phi}_{n}$ with respect to injectivity or
$\boldsymbol{\Phi}_{n}(B_{2})\subset\Omega$ a similar plot of
$\boldsymbol{\Phi}_{n}(B_{3})$ is not possible. For this reason, at the end of
each optimization we calculate two measures which help us to decide if the
constructed $\boldsymbol{\Phi}_{n}$ is injective and into.

On the other hand the principal approach to constructing $\boldsymbol{\Phi
}_{n}$ is very similar to the algorithm described in \S \ref{iteration}. Again
we are looking for a function $\boldsymbol{\Phi}_{n}$ given in the following
form
\[
\boldsymbol{\Phi}_{n}(x)=\sum_{j=1}^{N_{n}}\alpha_{n,j}\psi_{j}(x),\quad\quad
x\in B_{3},
\]
where $\{\psi_{1},\ldots,\psi_{N_{n}}\}$ is an orthonormal basis of $\Pi_{n}$
and the vectors $\alpha_{n,j}\in\mathbb{R}^{3}$, $j=1,\ldots,N_{n}$ are
determined by an optimization algorithm.

For a given $n\in\mathbb{N}$ we use the extremal points of Womersley, see
\cite{WOM}, on the sphere $S^{2}$. We will denote these points by
$W_{n}=\{z_{1}^{(n)},\ldots,z_{(n+1)^{2}}^{(n)}\}$. These points guarantee
that the smallest singular value of the interpolation matrix
\[
A_{n}:=:\left[
\begin{array}
[c]{ccc}%
\psi_{1}(z_{1}^{(n)}) & \ldots & \psi_{N_{n}}(z_{1}^{(n)})\\
\vdots &  & \vdots\\
\psi_{1}(z_{(n+1)^{2}}^{(n)}) & \ldots & \psi_{N_{n}}(z_{(n+1)^{2}}^{(n)})
\end{array}
\right]
\]
stays above $1$ for all $n$ which we have used for our numerical examples. The
number $(n+1)^{2}$ is also the largest possible number of interpolation points
on the sphere which we can use, because $\dim(\Pi_{n}|_{S^{2}})=(n+1)^{2}$,
see \cite[Corollary 2.20 and formula (2.9)]{atkinson-han2011}. Again we
enforce
\[
\boldsymbol{\Phi}_{n}(z_{j}^{(n)})=\varphi(z_{j}^{(n)}),\quad\quad
j=1,\ldots,(n+1)^{2},
\]
for the mapping function $\boldsymbol{\Phi}_{n}$; see also (\ref{F4}). To
define the initial function
\begin{equation}
\boldsymbol{\Phi}_{n}^{(0)}(x)=\sum_{j=1}^{N_{n}}\alpha_{n,j}^{(0)}\psi
_{j}(x),\quad\quad x\in B_{3},\label{eq4000}%
\end{equation}
we choose
\begin{equation}
\alpha_{n,j}^{(0)}=(\widetilde{\Phi},\psi_{j})_{B_{3}},\quad\quad
j=1\ldots,N_{n}.\label{eq4001}%
\end{equation}
$(\cdot,\cdot)_{B_{3}}$ is the usual $L_{2}$ inner product on $B_{3}$. The
polynomial $\boldsymbol{\Phi}_{n}^{(0)}$ is the orthogonal projection of
$\widetilde{\Phi}$ into $\Pi_{n}$. The function $\widetilde{\Phi}$ is some
continuous extension of $\varphi$ to $B_{3}$, obtained by the generalization
to three dimensions of one of the methods discussed in
\S \S \ref{constructions},\ref{integinterp}. Having determined
$\boldsymbol{\Phi}_{n}^{(0)}$, we convert the constrained optimization of the
objective function $\Lambda(\cdot)$ into an unconstrained minimization, as
discussed earlier in (\ref{F4})-(\ref{F7}).

Once the \textsc{Matlab} program \textit{fminunc} returns a local minimum for
$\Lambda\left(  \boldsymbol{\Phi}_{n}\right)  $ and an associated minimizer
$\boldsymbol{\Phi}_{n}$, we need to determine if $\boldsymbol{\Phi}_{n}$
satisfies
\begin{align}
\boldsymbol{\Phi}_{n}(x)  &  \neq\boldsymbol{\Phi}_{n}(y),\quad\quad x,y\in
B_{3},\quad x\neq y,\quad\mbox{(injective)}\label{eq5001}\\
\boldsymbol{\Phi}_{n}(B_{3})  &  \subset\Omega,\quad\quad
\mbox{(into)}\label{eq5002}%
\end{align}
For this reason we calculate two measures of our mapping $\boldsymbol{\Phi
}_{n}$.

Given $K\in\mathbb{N}$ we define a grid on the unit sphere,
\[%
\begin{array}
[c]{r}%
S_{K}^{2}:=\left\{  \left(  \sin\left(  \frac{1}{K}\pi j\right)  \cos\left(
\frac{1}{K}i\pi\right)  ,\sin\left(  \frac{1}{K}\pi j\right)  \sin\left(
\frac{1}{K}i\pi\right)  ,\cos\left(  \frac{1}{K}\pi j\right)  \right)
\,|\,\right.  \smallskip\\
\left.  j=0,\ldots,K,\quad i=0,\ldots,2K-1^{\,^{\left.  {}\right.  }}\right\}
.
\end{array}
\]
For $L\in\mathbb{N}$, we define a cubic grid in $B_{3},$%
\[
B_{3,L}:=\left(  \frac{1}{L}\mathbb{Z}^{3}\right)  \cap B_{3}%
\]
so every element in $B_{3,L}$ is given by
\[%
\begin{array}
[c]{c}%
\frac{1}{L}(i,j,k),\quad\quad i,j,k\in\mathbb{Z},\smallskip\\
i^{2}+j^{2}+k^{2}\leq L^{2}.
\end{array}
\]
To measure the minimum of the magnitude of the gradient of $\varphi$ over
$S^{2}$, we define an approximation by
\[
m_{K}(\varphi):=\displaystyle\min_{\substack{x,y\in S_{K}^{2} \\x\neq y}%
}\frac{\Vert\varphi(x)-\varphi(y)\Vert}{\Vert x-y\Vert}%
\]
This number is used to calculate
\[
E_{1,K}(\boldsymbol{\Phi}_{n}):=\displaystyle\min_{\substack{x,y\in B_{3,L},
\\x\neq y}}\frac{\Vert\boldsymbol{\Phi}_{n}(x)-\boldsymbol{\Phi}_{n}(y)\Vert
}{\Vert x-y\Vert}\Big/m_{K}(\varphi)
\]
Because of $\boldsymbol{\Phi}_{n}|_{S^{2}}\approx\varphi$ we expect
$E_{1,K}\leq1$. We use the occurrence of a very small value for $E_{1,K}%
(\boldsymbol{\Phi}_{n})$ to indicate that (\ref{eq5001}) may be violated. The
result $E_{1,K}(\boldsymbol{\Phi}_{n})\approx1$ is the best we can achieve,
for example, with $\varphi$ and $\boldsymbol{\Phi}_{n}$ the identity mapping.

If (\ref{eq5002}) is violated there is a point $x\in B_{3}$ and a point $y\in
S^{2}$ with $\boldsymbol{\Phi}_{n}(x)=\varphi(y)$. This shows that the
following measure would be close to zero
\[
E_{2,K,L}(\boldsymbol{\Phi}_{n}):=\displaystyle\min_{x\in B_{3,L},\,y\in
S_{K}^{2}}\frac{\Vert\boldsymbol{\Phi}_{n}(x)-\varphi(y)\Vert}{\Vert x-y\Vert
}\Big/m_{K}(\varphi)
\]
Again we expect $E_{2,K,L}(\boldsymbol{\Phi}_{n})\leq1$, and a very small
value of $E_{2,K,L}(\boldsymbol{\Phi}_{n})$ indicates that (\ref{eq5002})
might be violated. For each $\boldsymbol{\Phi}_{n}$ which we calculate we will
always report $E_{1,K}(\boldsymbol{\Phi}_{n})$ and $E_{2,K,L}(\boldsymbol{\Phi
}_{n})$. For larger $K$ and $L$ we will get a more accurate test of the
conditions (\ref{eq5001}) and (\ref{eq5002}), but the cost of calculation is
rising, the complexity to calculate $E_{2,K,L}(\boldsymbol{\Phi}_{n})$ for
example is $O(n^{3}K^{2}L^{3})$. For our numerical results we will use $K=40$
and $L=10$.

We consider only starlike examples for $\Omega$, with $\partial\Omega$ given
as%
\begin{align}
\varphi\left(  x\right)   &  =\rho\left(  x\right)  x,\quad\quad x\in
S^{2}\nonumber\\
&  =\widehat{\rho}\left(  \theta,\phi\right)  \left(  \sin\theta\cos\phi
,\sin\theta\sin\phi,\cos\theta\right) \label{eq5003}%
\end{align}%
\[
\widehat{\rho}\left(  \theta,\phi\right)  =\rho\left(  \sin\theta\cos\phi
,\sin\theta\sin\phi,\cos\theta\right)
\]
To create an initial guess, we begin with the generalization of (\ref{eq6}%
)-(\ref{eq7}) to three dimensions, defined in the following way:%
\begin{align*}
\widetilde{\Phi}\left(  x\right)   &  =r\widehat{\rho}\left(  \theta
,\phi\right)  \left(  \sin\theta\cos\phi,\sin\theta\sin\phi,\cos\theta\right)
\\
&  =\widehat{\rho}\left(  \theta,\phi\right)  x
\end{align*}
for $x=r\left(  \sin\theta\cos\phi,\sin\theta\sin\phi,\cos\theta\right)  $,
$0\leq\theta\leq\pi,\ 0\leq\phi\leq2\pi,\ 0\leq r\leq1$. We assume $\rho
:S^{2}\rightarrow\partial\Omega$ is a given smooth positive function. The
initial guess $\boldsymbol{\Phi}_{n}^{\left(  0\right)  }$ is obtained using
(\ref{eq4000})-(\ref{eq4001}), the orthogonal projection of $\widetilde{\Phi}$
into $\Pi_{n}$. Even though $\widetilde{\Phi}$ is not continuously
differentiable over $B_{3}$, its orthogonal projection $\boldsymbol{\Phi}%
_{n}^{\left(  0\right)  }$ is continuously differentiable, and it turns out to
be a suitable initial guess with $\left.  \boldsymbol{\Phi}_{n}^{\left(
0\right)  }\right\vert _{S^{2}}\approx\varphi$.

\begin{example}
\label{ex5001} In our first example we choose
\begin{equation}
\widehat{\rho}(\theta,\phi):=2+\left(  \cos\theta\right)  ^{2}\label{eq5005}%
\end{equation}
Using $n=6$ yields the results given in Table \ref{Table-3D-1} for the mapping
$\Phi_{6}$ obtained using the optimization procedure described above.%

\begin{table}[tbp] \centering
\caption{Measures of approximation stability for (\ref{eq5005})}\label{Table-3D-1}%
\begin{tabular}
[c]{|c|c|c|}\hline
$\Lambda(\Phi_{6})$ & $E_{1,40}(\Phi_{6})$ & $E_{2,40,10}(\Phi_{6})$\\\hline
$3.0575574308$ & $0.7485506872$ & $0.6626332145$\\\hline
\end{tabular}%
\end{table}%
%

\begin{figure}[tb]%
\centering
\includegraphics[
height=5.0004in,
width=6.0001in
]%
{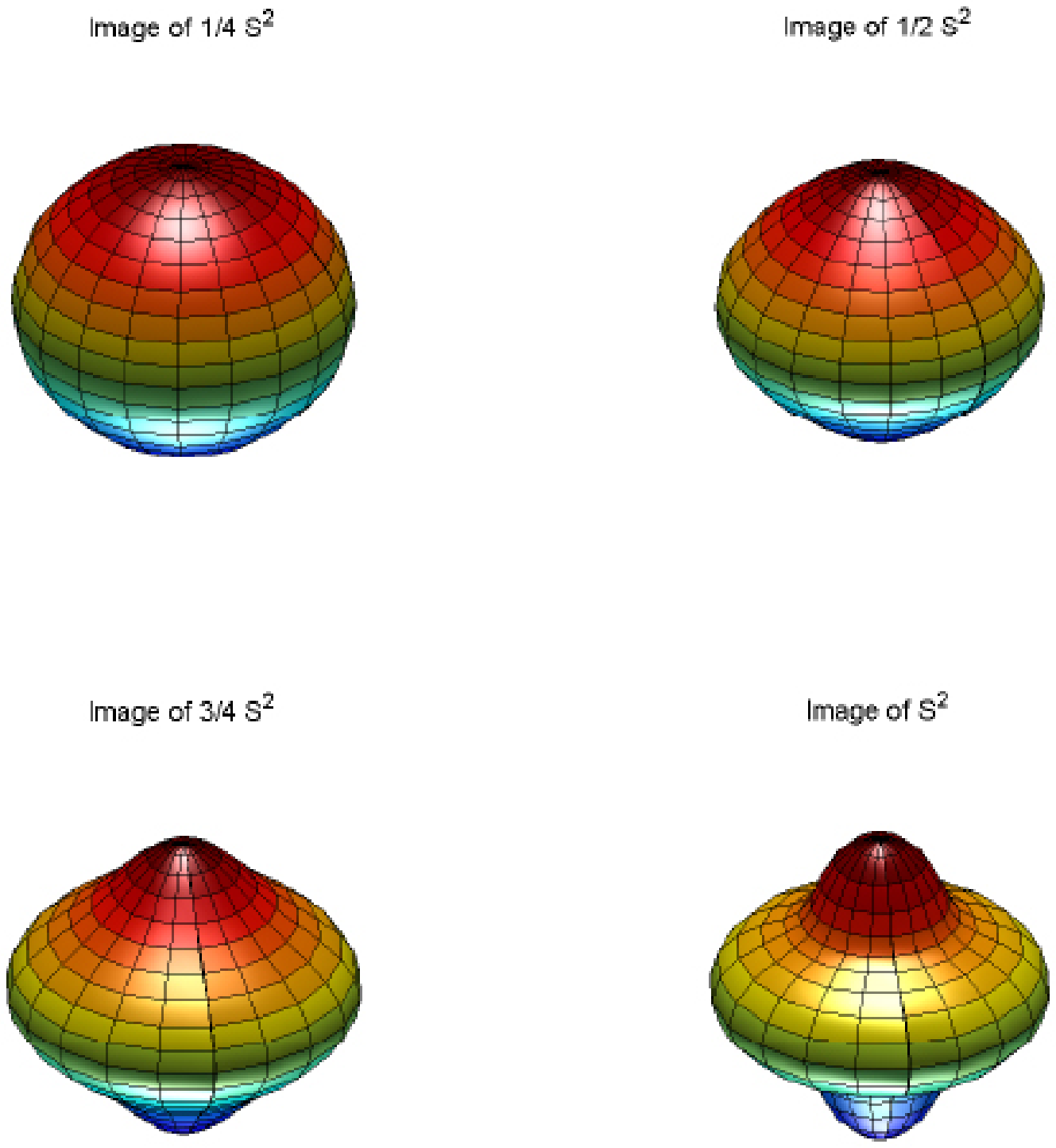}%
\caption{Images of $\frac{i}{4}S^{2}$, $i=1,2,3,4$}%
\label{fig5001}%
\end{figure}

See Figure \ref{fig5001} for an illustration of the images of the various
spheres $\frac{i}{4}S^{2}$. \ In this example the initial mapping
$\boldsymbol{\Phi}_{n}^{(0)}$ turned out to be a local optimum, so after the
first iteration the optimization stopped. The measures $E_{1}$ and $E_{2}$
seem to indicate that the function $\boldsymbol{\Phi}_{n}^{(0)}$ is into
$\Omega$ and injective. The error of $\Phi_{6}$ on the boundary is zero.
\end{example}

\begin{example}
\label{ex5002} Again the boundary $\partial\Omega$ is given by (\ref{eq5003}),
but this time we choose
\begin{equation}
\widehat{\rho}(\theta,\phi):=2+\cos\theta+\frac{1}{2}\sin\theta\sin
\phi\label{eq5006}%
\end{equation}
Using $n=6$ gives us the results shown in Table \ref{Table-3D-2}. We let
$\Phi_{6}^{(0)}$ denote our initial guess for the iteration, derived as
discussed earlier.%

\begin{table}[tbp] \centering
\caption{Measures of approximation stability for (\ref{eq5006})}\label{Table-3D-2}%
\begin{tabular}
[c]{|c|c|c|c|}\hline
Function & $\Lambda(\cdot)$ & $E_{1,40}(\cdot)$ & $E_{2,40,10}(\cdot)$\\\hline
$\Phi_{6}^{(0)}$ & $394.3717406299$ & $0.2088413520$ & $0.5926402745$\\\hline
$\Phi_{6}$ & $43.8782117161$ & $0.2018029407$ & $0.5175844592$\\\hline
\end{tabular}%
\end{table}%
%

\begin{figure}[tb]%
\centering
\includegraphics[
height=5.0004in,
width=6.0001in
]%
{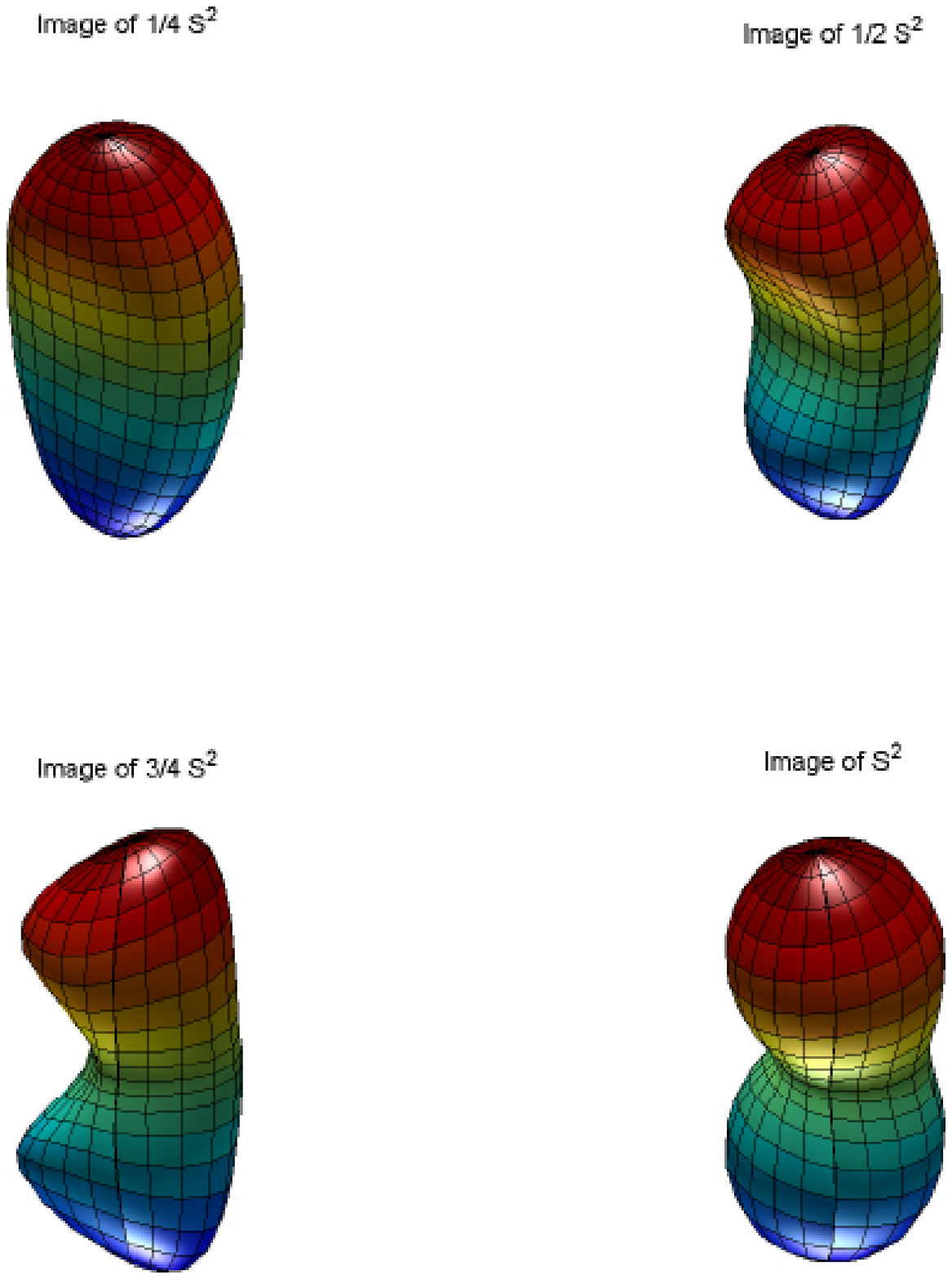}%
\caption{The images $\Phi^{\left(  34\right)  }\left(  \frac{i}{4}%
S^{2}\right)  ,$ $i=1,2,3,4$}%
\label{fig5003}%
\end{figure}

See Figure \ref{fig5003} for an illustration of the images of the various
spheres $\frac{i}{4}S^{2}$. In this example the $\Lambda(\cdot)$ value of the
initial mapping $\Phi_{6}^{(0)}$ is significantly improved by the
optimization. During the optimization the measures $E_{1}$ and $E_{2}$ do not
approach zero, which indicates that $\Phi_{6}$ is a mapping from $B_{3}$ into
$\Omega$ and is injective. The error of $\Phi_{6}^{(0)}$ and $\Phi_{6}$ on the
boundary is zero.
\end{example}


\begin{thebibliography}{99}                                                                                               %
\bibitem {ACH2010}K. Atkinson, D. Chien, and O. Hansen. A Spectral Method for
Elliptic Equations: The Dirichlet Problem, \textit{Advances in Computational
Mathematics, }\textbf{33} (2010), pp. 169-189, DOI=10.1007/s10444-009-9125-8.

\bibitem {atkinson-han2004}K. Atkinson and W. Han. On the numerical solution
of some semilinear elliptic problems, \textit{Electronic Transactions on
Numerical Analysis} \textbf{17} (2004), pp. 206-217.

\bibitem {atkinson-han2011}K. Atkinson and W. Han. \textit{Approximation on
the Unit Sphere}, to be published.

\bibitem {AH2010}K. Atkinson and O. Hansen. A Spectral Method for the
Eigenvalue Problem for Elliptic Equations, \textit{Electronic Transactions on
Numerical Analysis }\textbf{37} (2010), pp. 386-412.

\bibitem {AHC2011}K. Atkinson,O. Hansen, and D. Chien. A Spectral Method for
Elliptic Equations: The Neumann Problem, \textit{Advances in Computational
Mathematics, }\textbf{34} (2011), pp. 295-317, \textit{\ }DOI=10.1007/s10444-010-9154-3.

\bibitem {Loga}B. Logan and L. Shepp. Optimal reconstruction of a function
from its projections, \textit{Duke Math. J.} \textbf{42}, (1975), 645--659.

\bibitem {Nocedal}J. Nocedal and S.J. Wright \textit{Numerical Optimization,}
Springer--Verlag New York, Inc., 1999.

\bibitem {stroud}A. Stroud. \textit{Approximate Calculation of Multiple
Integrals,} Prentice-Hall, Inc., Englewood Cliffs, N.J., 1971.

\bibitem {WOM}R. Womersley. Extremal (Maximum Determinant) points on the
sphere $\mathbb{S}^{2}$, http://web.maths.unsw.edu.au/$\sim$rsw/Sphere/Extremal/New/extremal1.html

\bibitem {ZEI}E. Zeidler. \textit{Nonlinear Functional Analysis and its
Applications I}, Springer--Verlag, New York Inc, 1986.
\end{thebibliography}
\end{document}